\documentclass[11pt,reqno]{amsproc}
\usepackage[margin=1in]{geometry}
\usepackage{amsmath, amsthm, amssymb}
\usepackage[colorlinks=true, pdfstartview=FitV, linkcolor=blue,citecolor=blue, urlcolor=blue]{hyperref}
\usepackage[abbrev,lite,nobysame]{amsrefs}
\usepackage{times}
\usepackage[usenames,dvipsnames]{color}

\def\dist{{\rm dist}}

\def\R {\mathbb{R}}

\DeclareMathOperator*{\esup}{ess\,sup}

\newtheorem{proposition}{Proposition}[section]
\newtheorem{theorem}[proposition]{Theorem}
\newtheorem{corollary}[proposition]{Corollary}
\newtheorem{lemma}[proposition]{Lemma}
\theoremstyle{definition}
\newtheorem{definition}[proposition]{Definition}
\newtheorem{remark}[proposition]{Remark}

\newtheorem{problem}{Problem}
\numberwithin{equation}{section}

\title[On non-autonomously forced Burgers equation]{On non-autonomously forced Burgers equation with periodic and Dirichlet boundary conditions}

\author[P. Kalita]{Piotr Kalita}
\address{Faculty of Mathematics and Computer Science, Jagiellonian University, ul. \L{}ojasiewicza 6, 30-348 Krak\'{o}w, Poland}
\email{piotr.kalita@ii.uj.edu.pl}

\author[P. Zgliczy\'{n}ski]{Piotr Zgliczy\'{n}ski}
\address{Faculty of Mathematics and Computer Science, Jagiellonian University, ul. \L{}ojasiewicza 6, 30-348 Krak\'{o}w, Poland}
\email{piotr.zgliczynski@ii.uj.edu.pl}

\thanks{Work of P.K. and P.Z. was supported by National Science Center (NCN) of Poland under project No. UMO-2016/22/A/ST1/00077, work of P.K. was also partially supported by NCN of Poland under projects No. DEC-2017/25/B/ST1/00302 and  DEC-2017/01/X/ST1/00408.}

\begin{document}
	\begin{abstract}
		We study the non-autonomously forced Burgers equation
		$$
		u_t(x,t) + u(x,t)u_x(x,t) - u_{xx}(x,t) = f(x,t)
		$$
		on the space interval $(0,1)$ with two sets of the boundary conditions: the Dirichlet and periodic ones.
		For both situations we prove that there exists the unique $H^1$ bounded trajectory of this equation defined for all $t\in \R$. Moreover we demonstrate that
		this trajectory attracts all trajectories both in pullback and forward sense. We also prove that for the Dirichlet case this attraction is exponential.
	\end{abstract}
\maketitle

\section{Introduction}
The questions about the attractor structure for dynamical systems governed by dissipative evolutionary partial differential equations (PDEs) are usually difficult, in many cases open, and belong to the key problems that are being researched in PDEs community. We focus on one situation where, as it turns out from our results, the structure of such attractor can be described fully. Namely, we study the asymptotic behavior for the following Burgers equation
$$
u_t(x,t) + u(x,t)u_x(x,t) - u_{xx}(x,t) = f(x,t),
$$
where $x\in (0,1)$ and the forcing is assumed to be non-autonomous. This equation serves as the most basic model which allows to understand the interference between the linear viscous term $-u_{xx}$ and the quadratic nonlinearity $u u_x$. We supply the equation with two sets of boundary conditions: the Dirichlet ones
$$ u(0,t) = u(1,t) = 0,$$
and the periodic ones
$$
u(0,t) = u(1,t)\quad \textrm{and}\quad u_x(0,t) = u_x(1,t).
$$
Assuming that the forcing $f$ belongs to the space $L^\infty(\R;L^2(0,1))$ we prove that for both cases the equation has a unique global in time trajectory which is uniformly bounded in time in $H^1$ norm and that this trajectory attracts all weak solutions both forward in time and in the pullback sense.

Our study starts with the a priori energy estimates, which follow the arguments of, e.g., \cites{LiTiti, Temam_2}. We note that in  \cite{LiTiti} the energy estimates and results on the solution regularity are derived for the unforced case. These estimates, together with the energy equation method, cf., e.g., \cite{Ball}, allow us to obtain the existence of the non-autonomous counterpart of the global attractor, namely the pullback attractor.  This object is a non-autonomous set, which attracts for a given time $t$, all the trajectories emanating from the bounded sets of initial data taken at time instants converging to minus infinity. The approach by pullback attractors to deal with asymptotic behavior of non-autonomous problems governed by PDEs started more than 20 years ago \cites{ChebKloeSchma, CrauelFlandoli} and has since then been used to study many classes of dissipative non-autonomous PDEs, see \cites{Luk, Lan, KloeRas, balibrea} for some the recent development of the theory.  We stress  that the pullback attractor existence for the considered problems is standard and needs only the energy methods. We provide the proofs, however, in order to make the article self contained, and moreover the results are used in the second part of the paper where we prove the global asymptotic stability of the unique eternal solution.  Using the argument inspired by the work of Hill and S\"{u}li \cite{suli} which uses the weak version of the maximum principle we prove that
the pullback attractor consists, in fact, of a single eternal trajectory. For the Dirichlet conditions, using the appropriate comparison principle, see, e.g., \cite{Friedman}, we prove that the attraction is  exponential in time. For periodic conditions while we expect that this attraction is exponential, we leave the question of the attraction speed, for now, open. We only prove that the unique eternal trajectory attracts, in forward and pullback sense, all trajectories, without obtaining the speed of attraction.

The problem with time independent $f$ and with the Dirichlet condition has been studied in \cite{suli}. The authors there prove that there exists the unique solution of the stationary problem which attracts all solutions of the evolutionary problem as time goes to infinity. They prove this for the case of the multidimensional domain $\Omega$ and for more general nonlinear term $a(u) \cdot \nabla u$. We note that such extensions of our present work are possible and straightforward, we chose to follow the one-dimensional problem only to avoid the technical bootstrapping arguments which are required, in the periodic case, to get the sufficient smoothness for the strong maximum principle.
We remark that the paper \cite{suli} only deals with the autonomous problem, and the question of the asymptotic behavior for the case of the Dirichlet condition and non-autonomous forcing was, to our knowledge, open. We fill this gap. We also remark that we strengthen even the autonomous result of \cite{suli} where the  time-independent forcing was assumed to be H\"{o}lder continuous $f\in C^\alpha(\overline{\Omega})$ and only the solutions with the initial data $u_0\in C(\overline{\Omega})$ were proved to be attracted (in the present work, as we consider only the one-dimensional case, $u_0\in L^2(0,1)$). 

The result with periodic boundary conditions is due to Jauslin, Kreiss, and Moser \cite{JauKreMos} who prove that if the forcing term is time periodic then there exists the unique eternal trajectory which attracts in future all solutions of the problem. The work of
\cite{JauKreMos} was later extended in \cites{Fontes, koreans}, where, always, the time periodicity of the forcing term was assumed. We remove this  periodicity requirement and show that the unique eternal solution attracts all trajectories for arbitrary $L^\infty(\R;L^2)$ forcing. Note, that in \cites{koreans} the authors use the energy method only, and not the maximum principle, and get the attraction only under the smallness assumption on the forcing term. They also provide the numerical evidence that for large forcing the unique eternal solution is not attracting all trajectories anymore. We prove that this is not the case and the unique eternal bounded solution is actually globally asymptotically stable independent on the forcing magnitude.

The direct motivation of our work are recent articles of  Cyranka and  Zgliczy\'{n}ski \cites{Cyranka, ZgliczynskiCyranka}. In \cite{Cyranka} the author obtained the existence of the globally asymptotically stable solution of the autonomous problem with the periodic boundary conditions, thus providing the computer assisted proof of the counterpart of the result of \cite{suli} with the Dirichlet conditions replaced by the periodic ones. On the other hand, in
\cite{ZgliczynskiCyranka} the authors obtained the existence of globally attracting solution for periodic boundary conditions and a non-autonomous and not necessarily periodic in time forcing having a given form. We underline that the advantage of \cites{Cyranka, ZgliczynskiCyranka} over the results of the present paper and of \cites{JauKreMos, suli} is that computer assisted methods do not need the maximum principle and they allow to construct more concrete bounds for the obtained attracting eternal solutions. Moreover, in the case of the periodic boundary conditions in \cite{ZgliczynskiCyranka} the exponential convergence speed is obtained, while we prove only global asymptotic stability, leaving the question of convergence speed in the general case, for now, open.

Our future aim is to construct the computer assisted technique in order to constructively obtain, with some accuracy, the unique attracting trajectory for the considered problem. To this end, in contrast to \cites{Cyranka, ZgliczynskiCyranka}, where the Fourier basis is used, we plan to use the Finite Element Method (FEM). The approach based on FEM is better suited to deal with the problems with Dirichlet conditions as the construction of the orthogonal basis, in the case of arbitrary multidimensional domain is in itself a hard problem. While the rigorous proofs obtained by means of computer assisted techniques obtained by FEM will be the topic of our forthcoming paper, here we focus on what can be obtained purely analytically. 

We also mention, that while we study only the problems with the Dirichlet and periodic boundary conditions, it appears very interesting to understand the asymptotic behavior with the Neumann conditions. Although this problem is no longer dissipative, in the unforced case, in \cite{CaoTiti}, Cao and Titi, prove that every trajectory converges to a stationary one. The same result was also obtained by Byrnes et al. \cite{Byrnes}, who use the infinite dimensional version of the center manifold theorem.

It also appears interesting to us, to extend the results of the present paper to study the global asymptotic stability for the non-autonomously forced Burgers equation with the fractional viscous term 
$$
u_t(x,t) + u(x,t)u_x(x,t) + \left(-\frac{\partial^2}{\partial x^2}\right)^\alpha u(x,t) = f(x,t).
$$
We hypothesize that, at least with periodic boundary condition, the result on the convergence to the unique eternal solution holds for $\alpha \in [1/2,1)$ as the two ingredients: regularization effect of the evolution and the maximum principle remain valid in this case of the ''weakened'' damping \cites{cordoba, kiselev}. This result is also suggested by the fact that the kernel of the fractional Burgers operator behaves similar as the kernel of the fractional Laplacian itself, cf. \cite{jakubowski}.

We end the introduction with the brief overview of our article structure. In Section \ref{sec:apriori} we derive the key energy estimates that we need to study our problems. We also discuss the existence, uniqueness, and regularity of the solutions. Section \ref{sec:pull} is devoted to the summary of required facts from the pullback attractors theory, and the results on the pullback attractors existence for the problems under consideration. Finally, in Section \ref{sec:last} we prove the uniqueness of the eternal trajectories as well as the results on the forward convergence, and, in the Dirichlet case, its speed.

\section{Problem formulation, strong and weak solutions and relations between them}\label{sec:apriori}
Throughout the paper we will denote by $C$ a generic positive constant which can change from line to line, sometimes even in the same formula the letter $C$ can appear several times and denote different constants. We denote $\Omega = (0,1)$, the space domain of problems under consideration. We will use the shorthand notation for the spaces of functions defined on $\Omega$, that is we will write $L^2 = L^2(0,1)$, $H^1_0 = H^1_0(0,1)$, $H^{-1} = H^{-1}(0,1)$, the dual space to $H^1_0$, and so on. For a Banach space $V$ we will denote by $\mathcal{P}(V)$, $\mathcal{B}(V)$ the family of, respectively, nonempty, and nonempty and bounded sets in $V$. The scalar product and norm in $L^2$ will be denoted by $(\cdot,\cdot)$ and $\|\cdot\|$, respectively.  For spaces other than $L^2$ we will always use the subscript to denote the corresponding norms and duality pairings.
By $\dot{H}^k$ we will denote the closure in $H^k$ norm of the space of restrictions to the interval $(0,1)$ of $1$-periodic functions belonging to $C^\infty(\R)$ such that their mean on the interval $(0,1)$ vanishes. If we do not impose the vanishing of the mean we denote the corresponding spaces by $H^k_{per}$.  We will frequently use the Poincar\'{e} inequality
$$
c\|v\| \leq \|v_x\|,
$$
valid for $v\in H^1_0$ and for $v\in \dot{H}^1$ with $c = 1/\pi^2$. We will also use the following well known interpolation inequalities,
\begin{align}
& \|v\|_{L^\infty} \leq \|v\|^{1/2}\|v_x\|^{1/2} \quad \textrm{for}\quad v\in H^1_0\ \ \textrm{or}\ \ v\in \dot{H}^1,\label{eq:interpol}\\
& \|v_x\| \leq \|v\|^{1/2}\|v_{xx}\|^{1/2} \quad \textrm{for}\quad v\in H^1_0\cap H^2 \ \ \textrm{or}\ \ v\in \dot{H}^2.\label{eq:interpol2}
\end{align}

Let $f \in L^\infty(R; L^2)$. We will always assume that $f$ is defined on the whole time axis $\mathbb{R}$ even though sometimes we will consider problems defined only on the interval $(t_0,\infty)$.
We will deal with two problems: the non-autonomously forced Burgers equation first with the Dirichlet and then with periodic boundary conditions. We start from the analysis of
the problem with the Dirichlet conditions. The main part od this section is devoted to the derivation of the energy estimates, cf. \cite{LiTiti}, where such estimates are derived for the unforced problem.

\subsection{Problem with the Dirichlet conditions}
We define the weak and strong form of the problem with the Dirichlet conditions.
\begin{problem}\label{prblm:weak_dir}
	Let $t_0 \in \mathbb{R}$, $u_0 \in L^2$ and $f\in L^\infty(\R;L^2)$. Find $u\in L^2_{loc}(t_0,\infty;H^1_0)$ with $u_t\in L^2_{loc}(t_0,\infty;H^{-1})$ such that
	\begin{align}
	& \langle u_t, v \rangle_{H^{-1}\times H^1_0} + (uu_x,v) + (u_x,v_x) = (f(t),v)\quad \textrm{for every}\quad v\in H^1_0\quad \textrm{and a.e.}\quad t\in (t_0,\infty),\label{formulation:weak}\\
	& u(t_0) = u_0.
	\end{align}
\end{problem}

\begin{problem}\label{prblm:strong_dir}
	Let $t_0 \in \mathbb{R}$, $u_0 \in H^1_0$ and $f\in L^\infty(\R;L^2)$. Find $u\in L^2_{loc}(t_0,\infty;H^2\cap H^1_0)$ with $u_t\in L^2_{loc}(t_0,\infty;L^2)$ such that
	\begin{align}
	& u_t + uu_x - u_{xx} = f(x,t)\quad \textrm{for almost every}\quad (x,t) \in (0,1)\times (t_0,\infty),\label{formulation:strong}\\
	& u(t_0) = u_0.
	\end{align}
\end{problem}

The proof of the following existence and uniqueness result is standard and follows by the Galerkin method, and hence we omit it. We only provide the key a priori estimates, which will be useful in the following part of the paper.
\begin{theorem}\label{thm:ex_uni}
	Problems \ref{prblm:weak_dir} and \ref{prblm:strong_dir} have unique solutions.
\end{theorem}
\begin{proof}
	Taking $v=u(t)$ in \eqref{formulation:weak} we obtain
	\begin{equation}\label{eq_basic1}
	\frac{1}{2}\frac{d}{dt}\|u(t)\|^2 + (u(t)u_x(t),u(t)) + \|u_x(t)\|^2 = (f(t),u(t)).
		\end{equation}
	Note that for a smooth function $v$ defined on the interval $(0,1)$ such that $v(0) = v(1)$ there holds
	$$
	\int_0^1 v(x) v(x) v_x(x)\, dx = \frac{1}{3}\int_0^1 \frac{d}{dx}v^3(x)\, dx = \frac{v^3(1)-v^3(0)}{3} = 0,
	$$
	and the relation holds for every $v\in H^1$ such that $v(0) = v(1)$ by the density of smooth functions in that space. Using this equality in \eqref{eq_basic1} after simple transformations we deduce 
	\begin{equation}\label{eq:1}
		\frac{d}{dt}\|u(t)\|^2 + \|u_x(t)\|^2 \leq  C\|f\|_{L^\infty(\R;L^2)}^2 \quad \textrm{for a.e.}\quad  t\in (t_0,\infty).
	\end{equation}
	By the Poincar\'{e} inequality and the Growall lemma we deduce
\begin{equation}\label{est:energy1}
	\|u(t)\|^2 \leq \|u(t_0)\|^2e^{-C(t-t_0)} + C\|f\|_{L^\infty(\R;L^2)}^2(1-e^{-C(t-t_0)}).
	\end{equation}
	It is also clear that
		\begin{equation}\label{eq:2}
		\int_{t_1}^{t_2}\|u_x(t)\|^2\, dt \leq  \|u(t_1)\|^2 + C\|f\|^2_{L^\infty(\R;L^2)}(t_2-t_1) \quad \textrm{for every}\quad  t_0\leq t_1 < t_2.
		\end{equation}
	Multiplying \eqref{formulation:strong} by $-u_{xx}$ and integrating over interval $(0,1)$ yields
	\begin{equation}\label{eq:energy_eqn}
	\frac{1}{2}\frac{d}{dt}\|u_x\|^2 + \|u_{xx}\|^2 = -\int_0^1 f(t)u_{xx}\, dt + \int_{0}^1 uu_x u_{xx}\, dx.
	\end{equation}
	It follows that
	$$
	\frac{1}{2}\frac{d}{dt}\|u_x(t)\|^2 + \|u_{xx}(t)\|^2 \leq \|f\|_{L^\infty(\R;L^2)}\|u_{xx}(t)\| + \|u(t)\|_{L^\infty}\|u_x(t)\|\|u_{xx}(t)\|.
	$$
	After obvious transformations which use \eqref{eq:interpol} and \eqref{eq:interpol2} we obtain
		$$
	\frac{1}{2}\frac{d}{dt}\|u_x(t)\|^2 + \|u_{xx}(t)\|^2 \leq C\|f\|_{L^\infty(\R;L^2)}^2+\frac{1}{4}\|u_{xx}(t)\|^2 + \|u(t)\|^{5/4}\|u_{xx}(t)\|^{7/4}.
	$$
	We use the Young inequality with $\varepsilon$ which yields
\begin{equation}\label{eq:35}
\frac{d}{dt}\|u_x(t)\|^2 + \|u_{xx}(t)\|^2 \leq C\|f\|_{L^\infty(\R;L^2)}^2+ C\|u(t)\|^{10}.
\end{equation}	
Using \eqref{est:energy1} it follows that (note that the constant $C$ is allowed to depend on $f$ but not on the initial data)
\begin{equation}\label{eq:4}
\frac{d}{dt}\|u_x(t)\|^2 + \|u_{xx}(t)\|^2 \leq C(1+\|u(t_0)\|^{10}).
\end{equation}
Now, as $u_x(t)$ is mean free, we can use the Poincar\'{e} inequality to deduce that
\begin{equation}\label{eq:4.5}
\frac{d}{dt}\|u_x(t)\|^2 + c\|u_{x}(t)\|^2 \leq C(1+\|u(t_0)\|^{10}),
\end{equation}	
where $c$ is the Poincar\'e constant. Applying the Gronwall lemma yields
\begin{equation}\label{eq:5}
\|u_x(t)\|^2 \leq \|u_x(t_0)\|^2e^{C(t_0-t)} + C(1+\|u(t_0)\|^{10}).
\end{equation}
Coming back to \eqref{eq:4} it follows that
\begin{align}
&\int_{t_1}^{t_2}\|u_{xx}(t)\|^2\, dt \leq \|u_x(t_1)\|^2 + C(t_2-t_1)(1+\|u(t_0)\|^{10})\nonumber\\
&\qquad  \leq \|u_x(t_0)\|^2 + C(1+t_2-t_1)(1+\|u(t_0)\|^{10})\quad \textrm{for}\quad t_0\leq t_1 < t_2.\label{eq:55}
\end{align}
Finally,
$$
\int_{t_1}^{t_2}\|u_t(t)\|^2\, dt \leq C\left((t_2-t_1)\|f\|^2_{L^\infty(\R;L^2)} + \int_{t_1}^{t_2}\|u_{xx}(t)\|^2\, dt + \int_{t_1}^{t_2}\|u(t)\|_{L^\infty}^2\|u_x(t)\|_{L^2}^2\, dt \right),
$$
and
\begin{equation}\label{eq:56}
\int_{t_1}^{t_2}\|u_t(t)\|^2\, dt \leq C\left((t_2-t_1) + \int_{t_1}^{t_2}\|u_{xx}(t)\|^2\, dt + (t_2-t_1) \esup_{t\in [t_1,t_2]}\|u_x(t)\|^4_{L^2} \right),
\end{equation}
The required regularity follows.
\end{proof}

We also observe the simple corollary which follows from the definition of the weak and strong solutions and Theorem \ref{thm:ex_uni}.

\begin{corollary}\label{cor:22}
	Let $u$ be a weak solution with the initial data taken at time $t_0$. If $t_1>t_0$ then the function $u|_{[t_1,\infty)}$ is the weak solution with the initial data $u(t_1)$ taken at $t_1$.
	If, in turn, $u$ is a strong solution with the initial data taken at time $t_0$ then, if $t_1>t_0$, the function $u|_{[t_1,\infty)}$ is the strong solution with the initial data $u(t_1)$.
	Moreover, if $u_0 \in H^1_0$ is the initial data taken at time $t_0$ then both the strong and weak solution with this initial data coincide.
\end{corollary}

In the next result we obtain the Lipschitz continuity on bounded sets of the mapping that assigns to the initial data the value of the strong solution after some time.

\begin{lemma}\label{lem:cont}
	If $u_0, v_0 \in H^1_0$ are the initial data taken at time $t_0$, such that $\|u_x(t_0)\|, \|v_x(t_0)\| \leq R$ and $u, v$ are strong solutions with these initial data, then for every $\tau>0$ there exists a constant $D(\tau, R)  > 0$ such that
	$$
	\|u_x(t_0+\tau) - v_x(t_0+\tau)\| \leq D(\tau, R) 	\|u_x(t_0) - v_x(t_0)\|.
	$$
\end{lemma}
\begin{proof}
	
		Let $u_0, v_0 \in H^1_0$ and let $u,v$ be strong solutions corresponding to $u_0, v_0$ at time $t_0$, respectively. Denoting $w=u-v$ there holds the following equation
	$$
	w_t - w_{xx}  + u u_x - vv_x  = 0 \quad  \textrm{a.e.}\quad t>t_0, x\in (0,1).
	$$
	Testing this equation by $-w_{xx}$, we obtain
	$$
	\frac{1}{2}\frac{d}{dt}\|w_x\|^2 + \|w_{xx}\|^2 \leq |(u w_x,w_{xx})| + |(v_xw,w_{xx})| .
	$$
	It follows that
	$$
	\frac{1}{2}\frac{d}{dt}\|w_x\|^2 + \|w_{xx}\|^2 \leq \|u\|_{L^\infty}\|w_x\|\|w_{xx}\| + \|w\|_{L^\infty}\|v_x\|\|w_{xx}\|.
	$$
	whence
	$$
	\frac{1}{2}\frac{d}{dt}\|w_x\|^2 + \|w_{xx}\|^2 \leq C \left(\|u_x\|+\|v_x\|\right)\|w_x\|\|w_{xx}\|.
	$$
	It follows that
	$$
	\frac{d}{dt}\|w_x\|^2 \leq C (\|u_x\|^2+\|v_x\|^2) \|w_x\|^{2}.
	$$
	The assertion follows by the Gronwall lemma and the estimate \eqref{eq:5}.
\end{proof}

Finally we will prove that the weak solution becomes instantaneously the strong one.

\begin{lemma}\label{lem:24}
	Let $u$ be a weak solution with the initial data $u(t_0) = u_0\in L^2$ taken at time $t_0$. Then $u(t) \in H^1_0$ for every $t>t_0$. Moreover $u|_{[t_0+\varepsilon;\infty)}$ is a strong solution with the initial data $u(t_0+\varepsilon)$ taken at time  $t_0+\varepsilon$. Finally
	for every set $B \in \mathcal{B}(L^2)$ and every $\varepsilon > 0$ there exists a set $B_\varepsilon \in \mathcal{B}(H^1_0)$ such that if $u$ is a weak solution with the initial data $u_0 \in B$ taken at time $t_0$, then $u(t) \in B_\varepsilon$ for every $t\geq t_0 + \varepsilon$.
\end{lemma}
\begin{proof}
	Again, the estimate that we derive is only formal. The actual estimate should be derived by considering the Galerkin solutions in the spaces spanned by the eigenfunctions of $-u_{xx}$ operator with the strongly in $L^2$ converging initial data.
	
	Coming back to \eqref{eq:4.5}, by the Gronwall lemma it follows that
	$$
	\|u_x(t_2)\|^2 \leq \|u_x(t_1)\|^2 + C(1+\|u(t_1)\|^{10}) \quad \textrm{for every}\quad t_0 \leq t_1 < t_2.
	$$
	Using \eqref{est:energy1} we deduce
	$$
	\|u_x(t_2)\|^2 \leq \|u_x(t_1)\|^2 + C(1+\|u(t_0)\|^{10}) \quad \textrm{for every}\quad t_0 \leq t_1 < t_2.
	$$
	Now we choose any $\varepsilon>0$ and integrate the above inequality with respect to $t_1$ over the interval $(t_0,t_0+\epsilon)$. It follows that
	$$
	\varepsilon\|u_x(t_2)\|^2 \leq \int_{t_0}^{t_0+\varepsilon}\|u_x(t_1)\|^2\, dt_1 + C\varepsilon(1+\|u(t_0)\|^{10}) \quad \textrm{for every}\quad t_0 +\varepsilon \leq  t_2.
	$$
	We can use \eqref{eq:2} to deduce
	$$
	\varepsilon\|u_x(t_2)\|^2 \leq  \|u(t_0)\|^2 +  C\varepsilon(1+\|u(t_0)\|^{10}) \quad \textrm{for every}\quad t_2 \geq t_0 + \varepsilon.
	$$
	Hence
	\begin{equation}\label{est:eqn2}
	\|u_x(t_2)\|^2 \leq  \frac{1}{\varepsilon}\|u(t_0)\|^2 +  C(1+\|u(t_0)\|^{10}) \quad \textrm{for every}\quad t_2 \geq t_0 + \varepsilon,
	\end{equation}
	and the assertion follows.
\end{proof}
\subsection{Problem with the periodic conditions} \label{sec:periodic_estimates} We will now consider the Burgers equation with periodic conditions $u(0,t) = u(1,t)$ and $u_x(0,t) = u_x(1,t)$. We assume that $f\in L^\infty(\R;\dot{L}^2)$, that is $\|f(t)\|_{L^2}$ is uniformly bounded, $f$ is $1$-periodic, and mean free. We define the weak and strong solutions as follows.
\begin{problem}\label{prblm:weak_per}
	Let $t_0 \in \mathbb{R}$, $u_0 \in \dot{L}^2$ and $f\in L^\infty(\R;\dot{L}^2)$. Find $u\in L^2_{loc}(t_0,\infty;\dot{H}^1)$ with $u_t\in L^2_{loc}(t_0,\infty;\dot{H}^{-1})$ such that
	\begin{align}
	& \langle u_t, v \rangle_{\dot{H}^{-1}\times \dot{H}^1} + (uu_x,v) + (u_x,v_x) = (f(t),v)\quad \textrm{for every}\quad v\in \dot{H}^1\quad \textrm{and a.e.}\quad t\in (t_0,\infty),\label{formulation:weak_per}\\
	& u(t_0) = u_0.
	\end{align}
\end{problem}

\begin{problem}\label{prblm:strong_per}
	Let $t_0 \in \mathbb{R}$, $u_0 \in \dot{H}^1$ and $f\in L^\infty(\R;\dot{L}^2)$. Find $u\in L^2_{loc}(t_0,\infty;\dot{H}^2)$ with $u_t\in L^2_{loc}(t_0,\infty;\dot{L}^2)$ such that
	\begin{align}
	& u_t + uu_x - u_{xx} = f(x,t)\quad \textrm{for almost every}\quad (x,t) \in (0,1)\times (t_0,\infty),\label{formulation:strong_per}\\
	& u(t_0) = u_0.
	\end{align}
\end{problem}

\begin{remark}
	Note that it is sufficient to restrict to the mean-free $f$ and the mean-free solution $u$ in periodic case. Indeed, suppose that $f\in L^\infty(\R;\dot{L}^2)$ is not necessarily mean-free. Denote  $\alpha(t) = \int_0^1 f(x,t)\, dx$ and $\beta(t) = \int_0^1 u(x,t)\, dx$. Taking $v = 1$ in \eqref{formulation:weak_per} we obtain
$$
\frac{d}{dt} \int_0^1 u(x,t) \, dx = \int_0^1 f(x,t)\, dt,
$$
hence $\beta(t)$ can be found by solving the ODE
$$
\beta'(t) = \alpha(t)\quad \textrm{with the initial data}\quad \beta(t_0) = \int_0^1 u_0(x)\, dx.
$$
Denote
$$
\gamma(t) = \int_{t_0}^t \beta(s)\, ds = (t-t_0) \int_0^1 u_0(x)\, dx + \int_{t_0}^t \int_{t_0}^s \alpha(r) \, dr\, ds,
$$
and
\begin{equation}\label{eq:mean_free}
v(x,t) = u(x-\gamma(t),t) - \beta(t).
\end{equation}
The function $v$ is $1$-periodic, similar as $u$, and it is mean-free. Moreover it satisfies the equation
$$
v_t + v_x \gamma' + \beta' - v_{xx} + (v+\beta)v_x = f.
$$
But $\gamma'(t) = \beta(t)$ and $\beta'(t) = \alpha(t)$, hence
$$
v_t - v_{xx} + vv_x = f(t) - \int_0^1 f(x,t)\, dx.
$$
The last equation can be solved for $v$ and \eqref{eq:mean_free} can be then used to recover $u$, the solution for the non-mean-free case.
\end{remark}

Similar as in the case with the Dirichlet conditions, the existence and uniqueness of the weak and strong solutions are standard, and follow by the Galerkin method. Hence we omit the proof of the next theorem, restricting only to giving the a priori estimates which are analogous to the ones in the Dirichlet case and will be needed in the subsequent computations.
 \begin{theorem}\label{thm:ex_uni_per}
 	Problems \ref{prblm:weak_per} and \ref{prblm:strong_per} have unique solutions.
 \end{theorem}
 \begin{proof}
 	Exactly as in the case of the Dirichlet condition taking $v=u(t)$ in \eqref{formulation:weak_per} we obtain
 	\begin{equation}\label{eq:2_per}
 	\int_{t_1}^{t_2}\|u_x(t)\|^2\, dt \leq  \|u(t_1)\|^2 + C\|f\|^2_{L^\infty(\R;\dot{L}^2)}(t_2-t_1) \quad \textrm{for every}\quad  t_0\leq t_1 < t_2,
 	\end{equation}
 	and, by the the Poincar\'{e} inequality for mean free functions and by the Gronwall lemma,
 	\begin{equation}\label{est:energy1_per}
 	\|u(t)\|^2 \leq \|u(t_0)\|^2e^{-C(t-t_0)} + C\|f\|_{L^\infty(\R;\dot{L}^2)}^2(1-e^{-C(t-t_0)}),
 	\end{equation}
 	whence from \eqref{eq:2_per} we deduce
 	\begin{equation} \label{eq:3_per}
 	\int_{t_1}^{t_2}\|u_x(t)\|^2\, dt \leq  \|u(t_0)\|^2 +  C\|f\|^2_{L^\infty(\R;\dot{L}^2)}(1+t_2-t_1) \quad \textrm{for every}\quad  t_0\leq t_1 < t_2.
 	\end{equation}
 	To derive the second energy inequality we multiply \eqref{formulation:strong_per} by $-u_{xx}$ and integrate over $(0,1)$ which yields
 	\begin{equation}\label{eq:energy_eqn_per}
 	\frac{1}{2}\frac{d}{dt}\|u_x\|^2 + \|u_{xx}\|^2 = -\int_0^1 f(t)u_{xx}\, dt + \int_{0}^1 uu_x u_{xx}\, dx.
 	\end{equation}
 	Proceeding exactly the same as in the case of the Dirichlet conditions, which is possible, as $u_{xx}(t)$ is mean free and hence we can use the Poincar\'{e} inequality, we deduce that
 	\begin{equation}\label{eq:4.5_per}
 	\frac{d}{dt}\|u_x(t)\|^2 + C\|u_{x}(t)\|^2 \leq C(1+\|u(t)\|^{10}),
 	\end{equation}	
 	and the estimate \eqref{est:energy1_per} as well as the Gronwall lemma yield
 	\begin{equation}\label{eq:5_per}
 	\|u_x(t)\|^2 \leq \|u_x(t_0)\|^2e^{C(t_0-t)} + C(1+\|u(t_0)\|^{10}).
 	\end{equation}
 	Analogously as in the Dirichlet case we also have the estimates
 	\begin{align}
 	&\int_{t_1}^{t_2}\|u_{xx}(t)\|^2\, dt \leq \|u_x(t_1)\|^2 + C(t_2-t_1)(1+\|u(t_0)\|^{10})\nonumber\\
 	&\qquad  \leq \|u_x(t_0)\|^2 + C(1+t_2-t_1)(1+\|u(t_0)\|^{10})\quad \textrm{for}\quad t_0\leq t_1 < t_2,\label{eq:55_per}
 	\end{align}
	and
 	\begin{equation}\label{eq:56_per}
 	\int_{t_1}^{t_2}\|u_t(t)\|^2\, dt \leq C\left((t_2-t_1) + \int_{t_1}^{t_2}\|u_{xx}(t)\|^2\, dt + (t_2-t_1) \esup_{t\in [t_1,t_2]}\|u_x(t)\|^4 \right).
 	\end{equation}
 \end{proof}
Exactly as in the Dirichlet case we have the result analogous to Corollary \ref{cor:22}.

\begin{corollary}\label{cor:22_per}
	Let $u$ be a weak solution with the initial data $u_0\in \dot{L}^2$ taken at time $t_0$. If $t_1>t_0$ then the function $u|_{[t_1,\infty)}$ is the weak solution with the initial data $u(t_1)$.
	If, in turn, $u$ is a strong solution with the initial data $u_0\in \dot{H}^1$ taken at time $t_0$ then, if $t_1>t_0$, the function $u|_{[t_1,\infty)}$ is the strong solution with the initial data $u(t_1)$.
	Moreover, if $u_0 \in \dot{H}^1$ is the initial data taken at time $t_0$ then both the strong and weak solution with this initial data coincide.
\end{corollary}

Similar as in the Dirichlet case the mapping which assigns to the initial data the value of the solution after given time is Lipschitz on bounded sets in $\dot{H}^1$. We skip the proof as it is analogous to the proof of Lemma \ref{lem:cont_per}.

\begin{lemma}\label{lem:cont_per}
	If, $u_0, v_0 \in H^1_0$ are the initial data taken at time $t_0$ such that $\|u_x(t_0)\|, \|v_x(t_0)\| \leq R$, and $u, v$ are strong solutions with these initial data, then for every $\tau>0$ there exists a constant $D(\tau, R)  > 0$ such that
	$$
	\|u_x(t_0+\tau) - v_x(t_0+\tau)\| \leq D(\tau, R)  	\|u_x(t_0) - v_x(t_0)\|.
	$$
\end{lemma}

The next result is analogous to the Lemma \ref{lem:24} for the Dirichlet case and the proof follows the same lines, so we skip it.

\begin{lemma}\label{lem:24_per}
	Let $u$ be a weak solution with the initial data $u_0\in \dot{L}^2$ taken at time $t_0$. Then $u(t) \in \dot{H}^1$ for every $t>t_0$. Moreover $u|_{[t_0+\varepsilon;\infty)}$ is a strong solution with the initial data $u(t_0+\varepsilon)$ taken at time  $t_0+\varepsilon$. Finally
	there holds the estimate
	\begin{equation}\label{est:eqn2_per}
	\|u_x(t_1)\|^2 \leq  \frac{1}{\varepsilon}\|u(t_0)\|^2 +  C(1+\|u(t_0)\|^{10}) \quad \textrm{for every}\quad \varepsilon>0 \quad \textrm{and}\quad t_1 \geq t_0 + \varepsilon.
	\end{equation}
\end{lemma}

\section{Pullback attractors and their existence} \label{sec:pull}
\subsection{Pullback attractors: definition and the result on existence.} We begin this section with the definition of a process and a pullback attractor.
\begin{definition}
	Let $V$ be a Banach space. A family of mappings $\{ S(t,t_0) \}_{t\geq t_0}$such that $S(t,t_0):V\to V$ is called a process if $S(t,t)$ is an identity for every $t\in \mathbb{R}$ and $S(t,t_1)S(t_1,t_0) = S(t,t_0)$ for every $t_0 \leq t_1 \leq t$.
\end{definition}

If $V$ is a Banach space and $A, B \subset V$, then the Hausdorff semidistance between these two sets is denoted by
$$
\dist_V(A,B) = \sup_{a\in A}\inf_{b\in B} \|a-b\|_V.
$$
We will call the families of sets $A(t) \in \mathcal{P}(V)$ parameterized by time $t\in \R$ \textit{non-autonomous} sets and denote them $\mathbb{A} = \{ A(t)\}_{ t\in \R }$.

We recall the definition of a pullback attractor.

\begin{definition}\label{defn:pullback}
	The non-autonomous set $\mathbb{A} = \{ A(t) \}_{t\in \R} $ is a \textit{pullback attractor} of a process $\{ S(t,t_0) \}_{t\geq t_0}$ on the Banach space $V$ if\begin{itemize}
		\item for every $t\in \R$ the set $A(t) \subset V$ is nonempty and compact,
		\item for every $s\geq t$ there holds $S(t,s) A (s) = A(t)$, i.e., the family $\mathbb{A}$ is invariant,
		\item for every $B \in \mathcal{B} (V)$ there holds
		$$
		\lim_{t_0 \to -\infty} \dist_V (S(t,t_0)B,A(t)) = 0,
		$$
		i.e., the family $\mathbb{A}$ is \textit{pullback attracting},
		\item if the non-autonomous set  $\mathbb{C} = \{ C(t)  \}_{t\in \mathbb{R}}$ is such that $C(t)$ is nonempty and compact for every $t \in \mathbb{R}$ and  $\mathbb{C}$ is pullback attracting, then $A(t) \subset C(t)$ for every $t\in \mathbb{R}$.
	\end{itemize}
\end{definition}
\begin{remark}
	It is straightforward to check that if there exists $B_0 \in \mathcal{B}(V)$ such that $A(t) \in B_0$ for every $t\in \mathbb{R}$ then the last assertion (minimality) follows from the first three.
\end{remark}

We also define the so called bounded eternal (complete) solutions and kernel sections of the process $ \{ S(t,t_0)\} _{t\geq t_0}$.
\begin{definition}
	The function $u:\mathbb{R} \to V$ is a bounded eternal solution of $ \{ S(t,t_0)\} _{t\geq t_0}$ if $\|u(t)\|_V \leq C$ for every $t\in \R$ and $S(t,t_0) u(t_0) = u(t)$ for every $t_0\in \R$ and $t\geq t_0$.
\end{definition}

\begin{definition}
	The non-autonomous set $\mathbb{K} = \{ K(t) \}_{t\in \R}$ is called a kernel section of $ \{ S(t,t_0)\} _{t\geq t_0}$ if
	$$
	K(t) = \{  u(t)\, :\ u\ \ \textrm{is a bounded eternal solution of}\  \{ S(t,t_0)\} _{t\geq t_0} \}.
	$$
\end{definition}

The existence of the pullback attractor and its relation with kernel sections follows from the next theorem. The proof, in a more general, multivalued, setting can be found for example in \cite{coti_kalita}.
\begin{theorem}\label{thm:pubblack}
	Suppose that the process $\{ S(t,t_0) \}_{t\geq t_0}$ on $V$ is such that
	\begin{itemize}
		\item the mappings $S(t,t_0):V\to V$ are continuous for every $t\geq t_0$,
		\item the process $ \{ S(t,t_0) \}_{t\geq t_0}$  is \textit{pullback asymptotically compact}, that is, for every $t\in \mathbb{R}$, every bounded sequence $ \{ x_n \} \subset V$ and every $t_n \to -\infty$ the sequence $S(t,t_n) x_n$ is relatively compact,
		\item the process $ \{ S(t,t_0) \}_{t\geq t_0}$  is \textit{pullback dissipative}, that is, there exists a set $B_0 \in \mathcal{B}(V)$ such that for every $B\in \mathcal{B}(V)$ and $t\in \mathbb{R}$ there exists $t_0=t_0(t,B)$ such that for every $t_1\leq t_0$ there holds $S(t,t_1)B \subset B_0$.
	\end{itemize}
Then $\{ S(t,t_0) \}_{t\geq t_0}$ has a pullback attractor $\mathbb{A} = \{ A(t) \}_{t\in \R}$ such that $A(t) \subset B_0$ for every $t\in \mathbb{R}$. This attractor is given by
$$
A(t) = \bigcap_{s\leq t} \overline{\bigcup_{\tau\leq s} S(t,\tau) B_0}.
$$
Moreover,  the pullback attractor $\mathbb{A}$ coincides with the kernel section $\mathbb{K}$.
\end{theorem}

\subsection{Pullback attractors: a bi-space attractor.} In the case of the Burgers equation the pullback attractor will be compact in $H^1_0$ and it will attract in the norm of $H^1_0$ all sets which are bounded in $L^2$. We will hence use the Babin's and Vishik's formalism of bi-space attractors, see \cite{Babin, cholewa_dlotko}. We assume that $H, V$ are two Banach spaces such that $V\subset H$ with a continuous embedding. The following definition of the bi-space pullback attractor differs from definition of the pullback attractor by requiring that it attracts not only the sets which are bounded in $V$ but also the sets which are bounded in $H$.
\begin{definition}
	Suppose that the family $\{ S(t,t_0) \}_{t\geq t_0}$ of mappings $S(t,t_0):H\to H$ is a process on $H$ and suppose that $S(t,t_0)|_V$, that is $S(t,t_0)$ restricted to $V$, is a process on $V$. The non-autonomous set $\mathbb{A} = \{ A(t) \}_{t\in \R} $ is a \textit{pullback $(H,V)$ attractor} of the process $\{ S(t,t_0) \}_{t\geq t_0}$ if\begin{itemize}
		\item for every $t\in \R$ the set $A(t) \subset V$ is nonempty and compact in $V$,
		\item for every $s\geq t$ there holds $S(t,s) A (s) = A(t)$, i.e., the family $\mathbb{A}$ is invariant,
		\item for every $B \in \mathcal{B} (H)$ there holds
		$$
		\lim_{t_0 \to -\infty} \dist_V (S(t,t_0)B,A(t)) = 0,
		$$
		i.e., the family $\mathbb{A}$ is \textit{pullback attracts} in $V$ the sets which are bounded in $H$,
		\item if the non-autonomous set  $\mathbb{C} = \{ C(t)  \}_{t\in \mathbb{R}}$ is such that $C(t)$ is nonempty and compact in $V$ for every $t \in \mathbb{R}$ and  $\mathbb{C}$  pullback attracts bounded sets in $V$ (such as in Definition \ref{defn:pullback}), then $A(t) \subset C(t)$ for every $t\in \mathbb{R}$.
	\end{itemize}
\end{definition}

 We prove the following theorem.
\begin{theorem}\label{thm:pullback_bi}
Suppose that the family $\{ S(t,t_0) \}_{t\geq t_0}$ of mappings $S(t,t_0):H\to H$ is a process on $H$ and suppose that $S(t,t_0)|_V$, that is $S(t,t_0)$ restricted to $V$, is a process on $V$ such that
	\begin{itemize}
	\item the mappings $S(t,t_0)|_V:V\to V$ are continuous in $V$ for every $t\geq t_0$,
	\item the process $ \{ S(t,t_0)|_V \}_{t\geq t_0}$  is \textit{pullback asymptotically compact} in $V$,
	\item the process $ \{ S(t,t_0)|_V \}_{t\geq t_0}$  is \textit{pullback dissipative} in $V$,
	\item for every $B\in \mathcal{B}(H)$ and for every $\epsilon>0$ the set $S(t+ \epsilon, t)B$ is bounded in $V$.
\end{itemize}
Then there exists the pullback $(H,V)$ attractor  $\mathbb{A} = \{ A(t) \}_{t\geq t_0}$ which coincides with the pullback attractor of  $\{ S(t,t_0)|_V \}_{t\geq t_0}$.
\end{theorem}
\begin{proof}
	In view of Theorem \ref{thm:pubblack} it only suffices to prove that $\mathbb{A}$ \textit{pullback attracts} in $V$ the sets which are bounded in $H$. So let $B\in \mathcal{B}(H)$. Then
	$$
	\dist_V (S(t,t_0)B,A(t)) = \dist_V (S(t,t_0+1)S(t_0+1,t_0)B,A(t)).
	$$
	The set $S(t_0+1,t_0)B$ is bounded in $V$ with the bound independent on the choice of $t_0$, cf. Lemma \ref{lem:24}. Hence
	$$
	\lim_{t_0 \to -\infty} \dist_V (S(t,t_0)B,A(t)) = \lim_{t_0+1\to -\infty} \dist_V (S(t,t_0+1)S(t_0+1,t_0)B,A(t)) = 0,
	$$
	and the proof is complete.
\end{proof}
\subsection{Existence of $(L^2,H^1_0)$ pullback attractor for the problem with the Dirichlet conditions.}\label{sec:pullback_dirichlet} We come back to the study of the Burgers equation with the Dirichlet conditions. In view of Corollary \ref{cor:22} and Lemma \ref{lem:24} the map
\begin{equation}\label{process_weak}
S(t,t_0)u_0 = \{  u(t)\, :\ u\ \textrm{is a weak solution of Problem}\ \ref{prblm:weak_dir}\ \textrm{with initial data}\ u_0\in L^2\ \textrm{at}\ t_0\}
\end{equation}
is a process on $L^2$, and the following relation holds
\begin{equation}\label{process_strong}
S(t,t_0)|_{H^1_0} u_0 = \{  u(t)\, :\ u\ \textrm{is a strong solution of Problem}\ \ref{prblm:strong_dir}\ \textrm{with initial data}\ u_0\in H^1_0\ \textrm{at}\ t_0\}.
\end{equation}
According to Lemma \ref{lem:cont} mappings $S(t,t_0)|_{H^1_0} :{H^1_0}\to {H^1_0}$ are continuous. Lemma \ref{lem:24} implies that for any $B\in \mathcal{B}(L^2)$ and for any $\varepsilon>0$ the set $S(t_0+\varepsilon,t_0)B$ belongs to $\mathcal{B}(H^1_0)$.
To get the existence of the $(L^2,H^1_0)$ attractor it is sufficient to
obtain the asymptotic compactness and dissipativity of the process given by the strong solutions. We start from the proof of dissipativity.
\begin{lemma}
	The process $\{S(t,t_0)|_{H^1_0}\}_{t\geq t_0}$ is pullback dissipative in $H^1_0$.
\end{lemma}
\begin{proof}
	Using \eqref{eq:35} and \eqref{est:energy1} it follows that
\begin{equation}\label{eq:35_a}
\frac{d}{dt}\|u_x(t)\|^2 + \|u_{xx}(t)\|^2 \leq C(1+\|u(t_0)\|^{10}e^{-C(t-t_0)}).
\end{equation}
By the Poincar\'{e} inequality we obtain
\begin{equation}\label{eq:35_b}
\frac{d}{dt}\|u_x(t)\|^2 + C\|u_{x}(t)\|^2 \leq C(1+\|u(t_0)\|^{10}e^{-C(t-t_0)}).
\end{equation}
Remembering that $C$ may denote three different constants in the above formula, after simple calculations which use the Gronwall lemma we obtain
\begin{align}
& \|u_x(t)\|^2 \leq \|u_x(t_0)\|^2e^{-C(t-t_0)} + C(1+\|u(t_0)\|^{10}e^{-C(t-t_0)})\nonumber\\
& \qquad  \leq \|u_x(t_0)\|^2e^{-C(t-t_0)} + C(1+\|u_x(t_0)\|^{10}e^{-C(t-t_0)}),\label{eq:6}
\end{align}
and the required dissipativity follows.
\end{proof}
There are several techniques to prove the asymptotic compactness. One of them relies on the existence of an absorbing set in a space compactly embedded in $H^1_0$. Since this technique would require additional regularity of $f$, to avoid the extra assumptions on $f$, we choose to use the technique based on the energy equation, cf, e.g., \cite{Ball}  in the proof of the next lemma.
\begin{lemma}\label{lem_as_comp}
	The process $\{S(t,t_0)|_{H^1_0}\}_{t\geq t_0}$ is pullback asymptotically compact in $H^1_0$.
\end{lemma}
\begin{proof}
Choose	$t\in \mathbb{R}$, and the sequence $t_n\to -\infty$, and a bounded sequence $ \{ u_{0n} \} \subset H^1_0$. Let $u_n$ be a strong solution corresponding to the initial data $u_{0n}$ taken at $t_n$. The estimate \eqref{eq:6} implies that $u_n(t-1)$ is a sequence bounded in $H^1_0$. We should prove that the sequence $u_n(t)$ is relatively compact in $H^1_0$.
Estimates \eqref{eq:55} and \eqref{eq:56} imply that the sequence $u_n$
is bounded in $L^2(t-1,t+1;H^2)$ and $L^\infty(t-1,t+1;H^1_0)$ and $u_{nt}$ is bounded in $L^2(t-1,t+1;L^2)$. The Aubin--Lions lemma implies that there exists $u\in L^\infty(t-1,t_1;H^1_0)\cap L^2(t-1,t+1;H^2)$, such that, for a nonrenumbered subsequence
\begin{align}
& u_n \to u\quad \textrm{strongly in}\  L^2(t-1,t+1;H^1_0) \ \textrm{and weakly}-* \ \textrm{in}\ L^\infty(t-1,t+1;H^1_0)\label{conv:1}\\
& u_{nt}\to u_{t}\quad \textrm{weakly in}\ L^2(t-1,t+1;L^2),\label{conv:timeder}\\
& u_{nxx}\to u_{xx}\quad \textrm{weakly in}\ L^2(t-1,t+1;L^2),\label{conv:weak}
\end{align}
and the last weak convergence also holds in $L^2(t-1,s;L^2)$ for every $s\in (t-1,t+1)$. In particular we deduce that
$$
u_n(s) \to u(s)\quad \textrm{for a.e.}\ s\in (t-1,t+1)\ \textrm{strongly in}\ H^1_0.
$$
Since $u\in L^2(t-1,t+1;H^2)\cap L^\infty(t-1,t+1;H^1_0)$ and $u_t\in L^2(t-1,t+1;L^2)$ it follows that $u\in C([t-1,t+1];H^1_0)$.
Consider the integrals
\begin{align*}
& \int_{t-1}^s (u_n(r)u_{nx}(r),u_{nxx}(r)) \, dr\\
&\quad  = \int_{t-1}^s ((u_n(r)-u(r))u_{nx}(r),u_{nxx}(r)) \, dr\\
& \qquad \qquad + \int_{t-1}^s (u(r)(u_{nx}(r)-u_x(t),u_{nxx}(r)) \, dr + \int_{t-1}^s (u(r)u_{x}(r),u_{nxx}(r)) \, dr
\end{align*}
for $s\in (t-1,t+1)$. The first two terms converge to zero as $n\to \infty$ due to the estimates
\begin{align*}
&\left|\int_{t-1}^s ((u_n(r)-u(r))u_{nx}(r),u_{nxx}(r)) \, dr\right| \leq
\left|\int_{t-1}^s \|u_n(r)-u(r)\|_{L^\infty}\|u_{nx}(r)\|_{L^2}\|u_{nxx}(r)\|_{L^2} \, dr\right|\\
&\qquad  \leq \sup_{r\in [t-1,t+1]}\|u_{nx}(r)\|_{L^2} \|u_n-u\|_{L^2(t-1,t+1;H^1_0)} \|u_n\|_{L^2(t-1,t+1;H^2)} \to 0 \quad \textrm{as}\quad n\to \infty,
\end{align*}
and
\begin{align*}
&\left|\int_{t-1}^s (u(r)(u_{nx}(r)-u_x(t),u_{nxx}(r)) \, dr\right| \leq
\left|\int_{t-1}^s \|u(r)\|_{L^\infty}\|u_{nx}(r)-u_x(r)\|_{L^2}\|u_{nxx}(r)\|_{L^2} \, dr\right|\\
&\qquad  \leq \sup_{r\in [t-1,t+1]}\|u_{x}(r)\|_{L^2} \|u_n-u\|_{L^2(t-1,t+1;H^1_0)} \|u_n\|_{L^2(t-1,t+1;H^2)} \to 0 \quad \textrm{as}\quad n\to \infty.
\end{align*}
It is straightforward to check that $uu_x$ in $L^2(t-1,t+1;L^2)$, whence
it follows that, for a subsequence still denoted by the same index,
$$
\lim_{n \to \infty}\int_{t-1}^s (u_n(r)u_{nx}(r),u_{nxx}(r)) \, dr = \int_{t-1}^s (u_n(r)u_{nx}(r),u_{nxx}(r)) \, dr.
$$
Coming back to \eqref{eq:energy_eqn} we deduce that the following energy equation holds for every $s\in [t-1,t+1]$
\begin{align*}
& \frac{1}{2}\|u_{nx}(s)\|^2 + \int_{t-1}^s (f(r),u_{nxx}(r))\, dr- \int_{t-1}^s (u_n(r)u_{nx}(r), u_{nxx}(r))\, dr\\
&\qquad  = \frac{1}{2}\|u_{nx}(t-1)\|^2 - \int_{t-1}^s\|u_{xx}(r)\|^2\, dr.
\end{align*}
Denote
$$
V_n(s) = \frac{1}{2}\|u_{nx}(t-1)\|^2 - \int_{t-1}^s\|u_{xx}(r)\|^2\, dr,
$$
and
$$
V(s) = \frac{1}{2}\|u_x(s)\|^2 + \int_{t-1}^s (f(r),u_{xx}(r))\, dr- \int_{t-1}^s (u(r)u_{x}(r), u_{xx}(r))\, dr,
$$
for $s\in [t-1,t+1]$. The functions $V_n$ are nonincreasing on $[t-1,t+1]$ and, since
$$
V_n(s) = \frac{1}{2}\|u_{nx}(s)\|^2 + \int_{t-1}^s (f(r),u_{nxx}(r))\, dr- \int_{t-1}^s (u_n(r)u_{nx}(r), u_{nxx}(r))\, dr,
$$
then
$$
V_n(s) \to V(s)\  \textrm{as}\ n\to \infty\ \textrm{for a.e.}\ s\in (t-1,t+1).
$$
Let $p_m\nearrow t$ and $r_m\searrow t$ be sequences such that $V_n(p_m)\to V(p_m)$ and $V_n(r_m)\to V(r_m)$.
Then
$$
V_n(p_m) \geq V_n(t) \geq V_n(r_m).
$$
Passing with $n$ to infinity it follows that
$$
V(p_m) \geq \limsup_{n\to \infty }V_n(t) \geq \liminf_{n\to \infty }V_n(t) \geq V(r_m).
$$
Passing with $m$ to infinity, continuity of $V$ implies that $V_n(1)\to V(1)$, whence $\|u_{nx}(t)\| \to \|u_x(t)\|$. Since the previous convergences \eqref{conv:1}--\eqref{conv:weak} imply that $u_n(t)\to u(t)$ weakly in $H^1_0$, we deduce that $u_n(t)\to u(t)$ strongly in $H^1_0$ and the proof is complete.
\end{proof}

We are in position to apply Theorem \ref{thm:pullback_bi} (and the fact that any solution in $L^2_{loc}(t_0,\infty;H^2\cap H^1_0)$ with time derivative in $L^2_{loc}(t_0,\infty;L^2)$ is a continuous function of time with values in $H^1_0$) to deduce the following result
\begin{theorem}\label{thm:pullback}
	Let $f\in L^\infty(\R;L^2)$. There exists a non-autonomous set $\mathbb{A} = \{ A(t) \}_{t\in \R}$, the pullback attractor for the process \eqref{process_strong} governed by the strong solutions, and
	$(L^2,H^1_0)$ pullback attractor for the process \eqref{process_weak} governed by the weak solutions.
	This attractor is given by
	$$
	A(t) = \{ u(t)\, :\ u \ \textrm{is a bounded in}\ H^1_0\ \textrm{eternal strong solution} \}.
	$$
	In particular this means that there exists at least one eternal strong solution of Problem \ref{prblm:strong_dir}. Moreover, each bounded eternal strong solution $u:\R\to H^1_0$ in $A(t)$ belongs to $C_b(\R;H^1_0)$.
\end{theorem}

	In the subsequent sections we will show that in fact the eternal strong solution is unique and in consequence the set $A(t)$ is a singleton for every $t\in \mathbb{R}$.
	
\subsection{Existence of $(\dot{L}^2,\dot{H}^1)$ pullback attractor for the problem with periodic conditions.} The argument of this section follows the lines of the argument for the Dirichlet case, so we skip the proofs, which are analogous to the ones in Section \ref{sec:pullback_dirichlet}. Similar as in the Dirichlet case, in view of Corollary \ref{cor:22_per} and Lemma \ref{lem:24_per} the map
\begin{equation}\label{process_weak_per}
	S(t,t_0)u_0 = \{  u(t)\, :\ u\ \textrm{is a weak solution of Problem}\ \ref{prblm:weak_per}\ \textrm{with initial data}\ u_0\in \dot{L}^2\ \textrm{at}\ t_0\}
\end{equation}
is a process on $\dot{L}^2$, and the following relation holds
\begin{equation}\label{process_strong_per}
	S(t,t_0)|_{\dot{H}^1} u_0 = \{  u(t)\, :\ u\ \textrm{is a strong solution of Problem}\ \ref{prblm:strong_per}\ \textrm{with initial data}\ u_0\in \dot{H}^1\ \textrm{at}\ t_0\}.
\end{equation}
The proof of the next theorem step by step follows the lines of the proof of Theorem \ref{thm:pullback}.  
\begin{theorem}\label{thm:pullback_per}
	Let $f\in L^\infty(\R;\dot{L}^2)$. There exists a non-autonomous set $\mathbb{A} = \{ A(t) \}_{t\in \R}$, the pullback attractor for the process \eqref{process_strong_per} governed by the strong solutions of the periodic problem, and
	$(\dot{L}^2,\dot{H}^1)$ pullback attractor for the process \eqref{process_weak_per} governed by its weak solutions.
	This attractor is given by
	$$
	A(t) = \{ u(t)\, :\ u \ \textrm{is a bounded in}\ \dot{H}^1\ \textrm{eternal strong solution} \}.
	$$
	In particular this means that there exists at least one eternal strong solution of Problem \ref{prblm:strong_per}. Moreover, each bounded eternal strong solution $u:\R\to \dot{H}^1$ belongs to $C_b(\R;\dot{H}^1)$.
\end{theorem}

In Section \ref{sec:pullback_max_periodic} we will prove that the convergence to the pullback attractor is in fact exponential and that for each $t$ the set $A(t)$ is a singleton.

\section{Convergence to the unique bounded eternal trajectory.} \label{sec:last}

 \subsection{Dirichlet conditions.} \label{sec:pullback_max_dirichlet} The argument of this section is inspired by \cite{suli}. Note, however, that in \cite{suli} the authors deal with the strong solutions. Their key result on the convergence to equilibrium, cf. \cite[Theorem 3.2]{suli}, is based on the comparison principle applied to the linear adjoint problem. This analysis of the linear problem does not depend on the fact if the original problem is autonomous or non-autonomous. We generalize, however, \cite[Theorem 3.2]{suli} because we combine the comparison principle with the Stampacchia argument \cite{stampacchia} which is valid even for weak solutions. Thus, we can consider more general class of forcing, while in \cite{suli} the authors require that $f$ is H\"{o}lder continuous. We also obtain our global asymptotic stability results for wider class of solutions, namely we allow that $u_0\in L^2$ and our solutions are not necessarily classical, but weak.
 
 We start the proof by showing that the eternal solution bounded in $H^1_0$ must be unique. This fact is established in the following theorem.
 
\begin{theorem}\label{theorem:l1}
Let $f\in L^\infty(\R;L^2)$. The pullback attractor $\mathbb{A} = \{ A(t) \}_{t\in \R}$ established in Theorem 	\ref{thm:pullback} consists of a single eternal solution $A(t) = \{u(t)\}$ for every $t\in \R$. In other words there exists a unique eternal solution $u \in C_b(\R;H^1_0)$ such that for every bounded set $B\in \mathcal{B}(L^2)$ of initial data there holds
$$
\lim_{t_0 \to -\infty} \dist_{H^1_0}(S(t,t_0)B,\{ u(t) \}) = 0.
$$
\end{theorem}
\begin{proof}
	Let $u,v: \R\to H^1_0$ be two eternal solutions such that $\|u(t)\|_{H^1_0} \leq M$ and $\|v(t)\|_{H^1_0} \leq M$ for every $t\in \R$. Denote $a(t) = \frac{1}{2}(v(t) + v(t))$. Then $\|a(t)\|_{H^1_0}\leq M$ and, by \eqref{eq:interpol}, $\|a(t)\|_{L^\infty} \leq M$ for every $t\in \R$. 
	Denote and $w(t) = u(t)-v(t)$. Then  $\|w(t)\|_{H^1_0}\leq 2M$ for every $t\in \R$. Moreover, the function $w \in C_b(\R;H^1_0)$ satisfies the equation
\begin{equation}\label{eq:forward}
	w_t - w_{xx} + (a w)_x = 0,
\end{equation}
	as well as the Dirichlet conditions $w(0,t) = w(1,t) = 0$ for every $t\in \R$.
	
	We must prove that $u(t) = v(t)$, i.e., $w(t) = 0$ for every $t\in R$. To this end we fix $t\in \R$. We will consider the above equation on time intervals $(t_0,t)$ for $t_0 < t$. First we observe that as $w(t) \in H^1_0 \subset C([0,1])$, the function $w(t)$ is continuous. Define two open sets $A_+$ and $A_-$ by $A_+ = \{ x\in [0,1]\, :\ w(x,t) > 0 \}$, and $A_- = \{ x\in [0,1]\, :\ w(x,t) < 0 \}$ and the function $z_0:[0,1] \to \R$ by
	$$
	z_0(x) = \chi_{A_+}(x) - \chi_{A_-}(x).
	$$
	It is clear that $z_0 \in L^2$. Moreover
	$$
	\|w(t)\|_{L^1} = \int_0^1 z_0(x) w(x,t)\, dx.
	$$
	Now, consider the backwards problem
\begin{align}
& z_t(s) + z_{xx}(s) + a(s)z_x(s) = 0\quad \textrm{for}\quad (x,s)\in (0,1)\times (t_0,t),\label{eq:backwards}\\
& z(0,s) = z(1,s) = 0\quad \textrm{for}\quad s \in [t_0,t],\\
& z(t) = z_0.\label{eq:backwards3}
\end{align}
It is standard to prove that this problem has a unique weak solution $z \in L^2(t_0,t;H^1_0)$ with $z_t\in L^2(t_0,t;H^{-1})$. Testing the weak form of \eqref{eq:backwards} with $w(s)$ we obtain
\begin{equation}\label{eq:7}
\langle z_t(s), w(s) \rangle_{H^{-1} \times H^1_0} - (z_x(s),w_x(s)) + (a(s)z_x(s), w(s)) = 0\quad \textrm{for a.e.}\quad s\in (t_0,t).
\end{equation}
If, in turn, we test \eqref{eq:forward} with $z(s)$, we arrive at
$$
( w_t(s), z(s) )  + (z_x(s),w_x(s)) + ( (a(s) w(s))_x, z(s)) = 0\quad \textrm{for a.e.}\quad s\in (t_0,t).
$$
Integrating by parts, it follows that,
$$
( w_t(s), z(s) )  + (z_x(s),w_x(s))  - ( a(s) w(s), z_x(s)) = 0\quad \textrm{for a.e.}\quad s\in (t_0,t).
$$
Adding \eqref{eq:7} to the last equation it follows that
$$
( w_t(s), z(s) )  + \langle z_t(s), w(s) \rangle_{H^{-1} \times H^1_0} = 0\quad \textrm{for a.e.}\quad s\in (t_0,t),
$$
whereas
$$
\frac{d}{dt}( w(s), z(s) ) = 0\quad \textrm{for a.e.}\quad s\in (t_0,t).
$$
Integrating the above inequality over the interval $(t_0,t)$ it follows that
\begin{equation}\label{eq:8}
(w(t_0),z(t_0)) = (w(t),z(t)) = (w(t),z_0) = \|w(t)\|_{L^1}.
\end{equation}
Introducing the time $\tau = t-s$ the problem \eqref{eq:backwards}--\eqref{eq:backwards3} is equivalent to the following forward in time problem
\begin{align*}
&y_\tau(\tau) - y_{xx}(\tau) - a(t-\tau) y_x(\tau) = 0 \quad \textrm{for}\quad (x,\tau) \in (0,1)\times (0,t-t_0),\\
&y(0,s) = y(1,s) = 0\quad \textrm{for}\quad s \in [0,t-t_0],\\
& y(0) = z_0,
\end{align*}
namely, its solution $y$ is given by $y(\tau) = z(t-\tau)$ for $\tau \in [0,t-t_0]$. We will use the comparison principle
\cite[Theorem 6.1]{Friedman}, also see \cite[Theorem 3.2]{suli}.
Define
$$
\psi(x,\tau) = \frac{e^{2(M+1)} - e^{x(M+1)}}{e^{2(M+1)}-e^{M+1}}e^{-\frac{\tau}{e^{2(M+1)}-1}}.
$$
Observe that $\psi(1,\tau) > 0$ and $\psi(0,\tau) > 0$ for $\tau \in [0,t-t_0]$. Since $\psi$ is smooth and $y$ satisfies the homogeneous Dirichlet condition at $x=0$ and $x=1$ this means that $(y(\tau)-\psi(\tau))_+$ belongs to $L^2(0,t-t_0;H^1_0)$. We calculate
\begin{align*}
& \psi_\tau(\tau) - \psi_{xx}(\tau) - a(t-\tau)\psi_x(\tau) \\
& \quad =
\frac{e^{-\frac{1}{e^{2(M+1)}-1} \tau}}{e^{2(M+1)}-e^{M+1}} \left( -\frac{e^{2(M+1)} - e^{x(M+1)}}{e^{2(M+1)}-1}  +  (M+1 + a(x,t-\tau) )(M+1) e^{x(M+1)} \right)\\
& \quad \geq \frac{e^{-\frac{1}{e^{2(M+2)}-1} \tau}}{e^{2(M+2)}-e^{M+2}} M > 0.
\end{align*}
We continue the argument using the weak maximum principle, in spirit of the method of Stampacchia \cite{stampacchia}. We obtain
$$
(\psi_{\tau}(\tau),(y(\tau)-\psi(\tau))_+) - (\psi_{xx}(\tau),(y(\tau)-\psi(\tau))_+) - (a(t-\tau)\psi_x(\tau),(y(\tau)-\psi(\tau))_+) \geq 0,
$$
for almost every $\tau \in (0,t-t_0)$. Moreover
$$
\langle y_{\tau}(\tau),(y(\tau)-\psi(\tau))_+\rangle_{H^{-1} \times H^1_0} - (y_{xx}(\tau),(y(\tau)-\psi(\tau))_+) - (a(t-\tau)y_x(\tau),(y(\tau)-\psi(\tau))_+) = 0,
$$
for almost every $\tau \in (0,t-t_0)$. It follows that
\begin{align*}
& \langle y_{\tau}(\tau)-\psi_{\tau}(\tau),(y(\tau)-\psi(\tau))_+\rangle_{H^{-1} \times H^1_0} - (y_{xx}(\tau)-\psi_{xx}(\tau),(y(\tau)-\psi(\tau))_+)\\
& \qquad  - (a(t-\tau)(y_x(\tau)-\psi_x(\tau)),(y(\tau)-\psi(\tau))_+) \leq  0.
\end{align*}
We deduce that
$$
\frac{1}{2}\frac{d}{dt}\|(y(\tau)-\psi(\tau))_+\|^2 + \left\|\frac{\partial}{\partial x}(y(\tau)-\psi(\tau))_+\right\|^2 \leq M \left\|\frac{\partial}{\partial x}(y(\tau)-\psi(\tau))_+\right\|\|(y(\tau)-\psi(\tau))_+\|,
$$
whence
$$
\frac{d}{dt}\|(y(\tau)-\psi(\tau))_+\|^2 \leq \frac{M^2}{2}\|(y(\tau)-\psi(\tau))_+\|^2\quad \textrm{for almost every}\quad \tau\in (0,t-t_0).
$$
Now, the Gronwall lemma implies that
$$
\|(y(\tau)-\psi(\tau))_+\|^2 \leq e^{\frac{M^2}{2}}\|(z_0-\psi(0))_+\|^2 \quad \textrm{for  every}\quad \tau\in [0,t-t_0].
$$
It is easy to see that $z_0(x) \leq 1$ and $\psi(x,0) \geq 1$ for every $x\in [0,1]$, and hence $(z_0-\psi(0))_+ = 0$ for every $x\in [0,1]$, and, in consequence $(y(\tau)-\psi(\tau))_+ = 0$ for every $\tau\in [0,t-t_0]$.
This means that
$$
y(x,t-t_0) \leq \psi(x,t-t_0)\quad\textrm{for a.e.}\quad x\in (0,1).
$$
In a similar way, testing by $(y(\tau)+\psi(\tau))_-$ in place of $(y(\tau)-\psi(\tau))_+$ it follows that
$$
-\psi(x,t-t_0) \leq y(x,t-t_0) \quad\textrm{for a.e.}\quad x\in (0,1).
$$
We deduce that
$$
\|z(t_0)\|_{L^\infty} = \|y(t-t_0)\|_{L^\infty} \leq \max_{x\in [0,1]} \psi(x,t-t_0) = C e^{-C(t-t_0)}.
$$
Coming back to \eqref{eq:8} we observe that
$$
\|w(t)\|_{L^1} \leq \|w(t_0)\|_{L^1} \|z(t_0)\|_{L^\infty} \leq C e^{-C(t-t_0)} \|w(t_0)\|_{L^1} \leq C e^{-C(t-t_0)} \|w(t_0)\|_{H^1_0} \leq 2MCe^{-C(t-t_0)}.
$$
By taking $t_0$ sufficiently negative it follows that for every $\varepsilon > 0$ there holds  $\|w(t)\|_{L^1} \leq \varepsilon$ and hence it has to be $w(t) = 0$. The proof is complete.
\end{proof}
 As a special case, when the set $B \in \mathbb{B}(L^2)$ is a singleton we obtain the following result.
 \begin{theorem}
 	Let $f\in L^\infty(\R;L^2)$. For every $u_0\in L^2$ there holds
 	$$
 	\lim_{t_0 \to -\infty} \|S(t,t_0)u_0 - u(t)\|_{H^1_0} = 0.
 	$$
 \end{theorem}
\begin{remark}
	If $f$ is independent of time, then the problem becomes autonomous and the process is actually a semigroup $\{S(t)\}_{t\geq 0}$. In such a case the above result states that if $f\in L^2$ then there exists the unique $u\in H^1_0$, the solution to the stationary problem, such that for every $u_0\in L^2$
	$$
	\lim_{t\to\infty}\|S(t)u_0 - u\|_{H^1_0} = 0,
	$$
	where $S(t)u_0$ is the value of the weak solution at time $t$ with the initial data $u_0$ taken at time equal to zero. Observe that we have strengthened the result of Hill and S\"{u}li \cite{suli} who require that $f\in C^\alpha([0,1])$ and who consider only classical solutions. Note, however, that in \cite{suli} the authors consider the case where the domain is not the interval $[0,1]$ but a bounded and open set $\Omega \subset \R^d$, and their nonlinear term has the form $a(u) \cdot \nabla u$.
\end{remark}
We can extend  Theorem \ref{theorem:l1} to get following result
\begin{lemma}\label{remark:l2convergence}
	 If $u:[t_1,\infty) \to H^1_0$ and $v:[t_1,\infty) \to H^1_0$ are two strong solutions such that $\|u_x(t)\| \leq M$ and $\|v_x(t)\| \leq M$ for every $t\geq t_1$, then
	$$
	\| u(t) - v(t)\|\leq C e^{-C(t-t_1)},
	$$
	where the constant $C$ depends on $M$.
\end{lemma}
\begin{proof}
In the course of the proof of Theorem \ref{theorem:l1} we have shown that
$$
\| u(t) - v(t)\|_{L^1} \leq C e^{-C(t-t_1)}.
$$
Now, for $v\in H^1_0$ by interpolation we get
$$
\|v\| \leq \|v\|_{L^\infty}^{1/2}\|v\|_{L^1}^{1/2}\leq \|v\|^{1/4}\|v_x\|^{1/4}\|v\|_{L^1}^{1/2}.
$$
It follows that
$$
\|v\| \leq \|v\|_{L^1}^{2/3}\|v_x\|^{1/3},
$$
and the assertion is proved.
\end{proof} 

\begin{remark}
	It is clear that if the non-autonomous forcing $f$ is $T$-periodic, then the pullback attractor $\mathbb{A} = \{ A(t) \}_{t\in \R}$ is also $T$-periodic, i.e. $A(t) = A(t+T)$ for every $t\in \R$, and, in our case, the unique eternal solution bounded in $H^1_0$ is periodic.
\end{remark}
For a set $B\in \mathcal{B}(L^2)$ we will denote by $\|B\|$ the value $\sup_{b\in B}\|b\|$. We prove the following result.
\begin{theorem}
	Let $f\in L^\infty(\R;L^2)$. Let $v_0 \in B \in \mathcal{B}(L^2)$ and let $v$ be a weak solution starting from the initial data $v_0$ at time $t_0$. Let $u$ be the unique eternal solution bounded in $H^1_0$. There exists a constant $C>0$, a constant $C(\|B\|)>0$ (depending continuously and monotonically on $\|B\|$) such that for every $t\geq t_0$ there holds
\begin{equation} \label{eqn:attraction}
	\|u(t) - v(t)\| \leq C(\|B\|) e^{-C(t-t_0)}.
\end{equation}
	In consequence, if only $f\in L^\infty(\R;L^2)$, then the unique eternal solution $u$ bounded in $H^1_0$ attracts exponentially in $L^2$ both in forward and pullback sense all weak solutions uniformly with respect to bounded sets of initial data in $L^2$.
\end{theorem}
\begin{proof}
	Let $v_0 \in L^2$ and let $v$ be a weak solution starting from $v_0$ at time $t_0$. Without loss of generality we consider only the case $\|v_0\| > 1$. If $\|v_0\|\leq 1$ we can take any fixed value greater than one in place of $\|v_0\|$. Assume that  $u:\R \to H^1_0$ is a unique eternal solution which is bounded in $H^1_0$. Estimate \eqref{est:energy1} implies that there exist positive constants $C, D$ such that if
	$$
	t \geq t_0 + C \ln \|v_0\|,
	$$
	then
	$$
	\|v(t)\| \leq D.
	$$
	Estimate \eqref{est:eqn2} with $\varepsilon = 1$ implies that there exist constants $C, D >0$ such that
	if
	$$
	t\geq 1+ t_0 + C \ln \|v_0\|,
	$$
	then
	$$
	\|v_x(t)\|\leq D.
	$$
	Consider $v|_{[1+ t_0 + C \ln \|v_0\|,\infty)}$ and  $u|_{[1+ t_0 + C \ln \|v_0\|,\infty)}$. Lemma \ref{remark:l2convergence} implies that for every $t\geq 1+ t_0 + C \ln \|v_0\|$ there holds
	$$
	\|u(t) - v(t)\| \leq Ce^{-C(t-1-t_0-\ln \|v_0\|)}.
	$$
	Hence, if only $t-t_0 \geq 1+ C \ln \|v_0\|$, then
	$$
	\|u(t) - v(t)\| \leq C\|v_0\|^C e^{-C(t-t_0)},
	$$
	and the assertion \eqref{eqn:attraction} follows for every $t\geq t_0 + 1+ C\ln \|v_0\|$.
	The fact that the assertion \eqref{eqn:attraction} holds for every $t\geq t_0$ follows from the fact that  $\|u(t)-v(t)\|$ is uniformly bounded in $L^2$.
\end{proof}

 \subsection{Periodic conditions.} \label{sec:pullback_max_periodic} Contrary to the Dirichlet case, in periodic situation we only prove that all trajectories converge (forward in $L^2$ and pullback in $H^1$) to a unique eternal strong solution bounded in $\dot{H}^1$. For the time, we leave open the question of the convergence speed, which we expect to be exponential. In the course of the proof, contrary to the Dirichlet case, the weak maximum principle appears to be insufficient to get the corresponding result, and we need to apply its strong version. This requires us to do the additional bootstrapping to get the desired regularity for the adjoint problem. Note that the convergence of higher order space derivatives is easy to obtain from our results, by interpolation and uniform a priori estimates in higher order norms which will hold under increased regularity assumptions on the forcing term $f$. 
 
 Before we pass to the proof of global asymptotic stability, we need two auxiliary results.

 \begin{lemma}\label{lemma:weak_max}
 	Assume that $a\in C_b([0,\infty);\dot{H}^1)$ satisfies $\|a(t)\|_{\dot{H}^1} \leq M$ for every $t\geq 0$ and $y \in L^2_{loc}(0,\infty;H^1_{per}) \cap C([0,\infty);{L}^2)$ with $y_t \in L^2_{loc}(0,\infty;(H^1_{per})')$ is the unique weak solution of the linear problem
 		\begin{align}
 	&y_\tau(\tau) - y_{xx}(\tau) - a(\tau) y_x(\tau) = 0 \quad \textrm{for}\quad (x,\tau) \in (0,1)\times (0,\infty),\label{eq:adjoint_1_abs}\\
 	&y(0,\tau) = y(1,\tau)\quad \textrm{and}\quad  y_x(1,\tau) = y_x(0,\tau) \quad \textrm{for}\quad \tau \in (0,\infty),\\
 	& y(0) = y_0,\label{eq:adjoint_3_abs}
 	\end{align}
 	with $y_0 \in L^\infty$. Then
 	$$
 	\|y(\tau)\|_{L^\infty} \leq \|y_0\|_{L^\infty} \quad \textrm{for every}\quad \tau\geq 0.
 	$$
 \end{lemma}
 \begin{proof}
 The proof follows by the method of Stampacchia \cite{stampacchia}, as in the Dirichlet case. Weak form of the considered equation is the following
 \begin{equation}\label{eq:adjoint_per}
 \langle y_\tau(\tau),z\rangle_{(H^{1}_{per})'\times H^1_{per}} + (y_x(\tau),z_x) - (a(\tau) y_x(\tau),z) = 0 \quad \textrm{for every}\quad z\in H^1_{per}\quad \textrm{a.e.}\quad \tau>0,
 \end{equation}with the initial data $y(0) = y_0 \in L^2$.
 It follows from the standard argument that this equation has a unique weak solution with the regularity given in the statement of the lemma.	
 We derive the maximum principle estimate for this equation. To this end, first test the above equation by $(y(\tau)-\|y_0\|_{L^\infty})_+$. This leads to the bound
 \begin{equation}\label{eq:Stampa}
 \frac{1}{2}\frac{d}{dt}\|(y(\tau)-\|y_0\|_{L^\infty})_+\|^2 + \|((y(\tau)-\|y_0\|_{L^\infty})_+)_x\|^2 \leq M \|((y(\tau)-\|y_0\|_{L^\infty})_+)_x\|\|(y(\tau)-\|y_0\|_{L^\infty})_+\|,
 \end{equation}
 and, since $\|(y(0)-\|y_0\|_{L^\infty})_+\| = 0$, the argument based on the Gronwall lemma implies that $\|(y(\tau)-\|y_0\|_{L^\infty})_+\| = 0$ for every $\tau \in [0,\infty)$. A similar argument based on testing the equation by $(y(\tau)+\|y_0\|_{L^\infty})_-$ leads to the conclusion that $\|(y(\tau)+\|y_0\|_{L^\infty})_-\| = 0$ and  in fact $\|y(\tau)\|_{L^\infty} \leq \|y_0\|_{L^\infty}$ for every $\tau \in [0,\infty)$.
 \end{proof}
We have proved that $L^\infty$ norm of the solution cannot exceed the $L^\infty$ norm of the initial data. In the next result we show, using the strong version of the maximum principle, that the $L^\infty$ inequality in the assertion of the last lemma must be in fact strict.
 \begin{lemma}\label{lem:strong_max}
		Assume that $a\in C_b([0,\infty);\dot{H}^1)$ with $a_\tau \in L^2_{loc}(0,\infty;L^2)$ satisfies $\|a(t)\|_{\dot{H}^1} \leq M$ for every $t\geq 0$ and $y \in L^2_{loc}(0,\infty;H^1_{per}) \cap C([0,\infty);{L}^2)$ with $y_t \in L^2_{loc}(0,\infty;(H^1_{per})')$ is the unique weak solution of the linear problem \eqref{eq:adjoint_1_abs}--\eqref{eq:adjoint_3_abs}
	with $y_0 \in L^\infty$. Assume that $y_0$ is not equal to a constant function for a.e. $x\in (0,1)$. Then
	$$
	\|y(\tau)\|_{L^\infty} <  \|y_0\|_{L^\infty} \quad \textrm{for every}\quad \tau\geq 0.
	$$
\end{lemma}
\begin{proof}
	We first need to establish the regularity of $y$ in order to use the strong maximum principle. Fix $T>0$. Testing \eqref{eq:adjoint_1_abs} by $-y_{xx}(\tau)$ we obtain
	\begin{equation}\label{eq:secondenergy_adjoint}
	\frac{1}{2}\frac{d}{dt}\|y_x(\tau)\|^2 + \|y_{xx}(\tau)\|^2 \leq M \|y_x(\tau)\|\|y_{xx}(\tau)\|\quad \textrm{for a.e.}\quad  \tau>0.
	\end{equation}It follows that
	$$
	\frac{d}{dt}\|y_x(\tau)\|^2  \leq \frac{1}{2}M^2 \|y_x(\tau)\|^2,
	$$
	whence for every $0 < \tau_1\leq \tau_2$
	$$
	\|y_x(\tau_2)\|^2 \leq e^{1/2 M^2 (\tau_2-\tau_1)} \|y_x(\tau_1)\|^2  \leq  e^{1/2 M^2 \tau_2} \|y_x(\tau_1)\|^2.
	$$
	We can integrate the above estimate with respect to $\tau_1$ over the interval $(0,\tau_2)$, which yields
	$$
	\tau_2 \|y_x(\tau_2)\|^2 \leq e^{1/2 M^2 \tau_2} \int_0^{\tau_2} \|y_x(\tau_1)\|^2\, d\tau_1,
	$$
	whence
	$$
	\|y_x(\tau_2)\|^2 \leq \frac{e^{1/2 M^2 \tau_2}}{\tau_2} \int_0^{\tau_2} \|y_x(\tau_1)\|^2\, d\tau_1  =    \frac{e^{1/2 M^2 \tau_2}}{\tau_2} \|y\|^2_{L^2(0,\tau_2;H^1_{per})}.
	$$
	Coming back to \eqref{eq:secondenergy_adjoint} it follows that
	$$
	\int_{\tau_2}^{\tau_3}\|y_{xx}(\tau)\|^2\, d\tau \leq M\int_{\tau_2}^{\tau_3} \|y_x(\tau)\|^2\, d\tau + \|y_x(\tau_2)\|^2 \quad \textrm{for every}\quad 0< \tau_2\leq \tau_3.
	$$
	We deduce that for every $\varepsilon \in (0,T)$ there holds
	$$
	y_x \in L^{\infty}(\varepsilon,T;L^2)\quad \textrm{and}\quad y_{xx}\in L^2(\varepsilon,T;L^2).
	$$
	Moreover, as $|a(x,t)|\leq M$ for every $(x,t)\in [0,1]\times [0,T]$ it holds that
	$$
	y_t\in L^2(\varepsilon,T;L^2).
	$$ We deduce that
	\begin{equation}\label{eq:Y_C}
	y\in C([\varepsilon,T];H^1_{per}),
	\end{equation}
	whence $y\in C([0,1]\times [\varepsilon,T])$ and $y(t)$ is $1$-periodic with respect to variable $x$ for every $t>0$.
	Now we differentiate \eqref{eq:adjoint_1_abs} with respect to $\tau$ and denote $p = y_{\tau}$. This function satisfies the equation
	\begin{equation}\label{eq:eq16}
	p_\tau(\tau) - p_{xx}(\tau) - a_\tau(\tau) y_x(\tau) - a(\tau) p_x(\tau) = 0,
	\end{equation}
	with the periodic boundary conditions. We first test this equation with $p(\tau)$. We obtain
	\begin{equation*}
	\frac{1}{2}\frac{d}{dt}\|p(\tau)\|^2 + \|p_x(\tau)\|^2 \leq \int_0^1 |a_\tau(x,\tau) p(x,\tau) y_x(x,\tau)|\, dx + M \|p_x(\tau)\| \|p(\tau)\|.
	\end{equation*}
	It follows that
	\begin{equation*}
	\frac{1}{2}\frac{d}{dt}\|p(\tau)\|^2 + \|p_x(\tau)\|^2 \leq  \|a_\tau(\tau)\| \|p(\tau)\|_{L^\infty} \|y_x(\tau)\| + M^2\|p(\tau)\|^2 + \frac{1}{4} \|p_x(\tau)\|^2.
	\end{equation*}
	We deduce
	\begin{equation}\label{eq:est17}
	\frac{d}{dt}\|p(\tau)\|^2 + \|p_x(\tau)\|^2 \leq  2\|a_\tau(\tau)\|^2 \|y_x(\tau)\|^2 + 2M^2\|p(\tau)\|^2.
	\end{equation}
	We integrate this inequality from $\tau_1>0$ to $\tau_2\in (\tau_1,T]$, whence
	$$
	\|p(\tau_2)\|^2\leq \|p(\tau_1)\|^2 + 2\int_{0}^T\|a_\tau(\tau)\|^2\, d\tau\  \esup_{\tau\in (\tau_1,T]}\|y_x(\tau)\|^2 + 2M^2\int_{\tau_1}^T\|p(\tau)\|^2\, d\tau.
	$$
	We choose $\tau_2 \geq 2\varepsilon$, and integrate the above inequality with respect to $\tau_1$ over the interval $(\varepsilon,2\varepsilon)$ whence
	$$
	\varepsilon\|p(\tau_2)\|^2\leq \int_{\varepsilon}^T\|p(\tau_1)\|^2\, d\tau_1 + 2 \varepsilon \int_{0}^T\|a_\tau(\tau)\|^2\, d\tau\  \esup_{\tau\in (\varepsilon,T]}\|y_x(\tau)\|^2 + 2M^2\varepsilon\int_{\varepsilon}^T\|p(\tau)\|^2\, d\tau.
	$$
	It follows that $y_\tau = p \in L^\infty(2\varepsilon,T;L^2)$. As
	$$
	\|y_{xx}(\tau)\| \leq \|y_\tau(\tau)\| + M \|y_x(\tau)\| = \|p(\tau)\| + M \|y_x(\tau)\|,
	$$
	it follows that $y_{xx} \in L^\infty(2\varepsilon,T;L^2)$ and hence $y_x\in L^\infty(2\varepsilon,T;L^\infty)$.
	Coming back to \eqref{eq:est17} we deduce that $y_{\tau x} = p_x \in L^2(2\varepsilon,T;L^2)$. Equation \eqref{eq:eq16} implies that $p_\tau = y_{\tau\tau} \in L^2(2\varepsilon,T;L^2)$.
	Finally we test \eqref{eq:eq16} with $-p_{xx}(\tau)$ which yields
	$$
	\frac{1}{2}\frac{d}{dt}\|p_x(\tau)\|^2 + \|p_{xx}(\tau)\|^2 \leq M \|p_x(\tau)\|\|p_{xx}(\tau)\| + \|a_\tau(\tau)\| \|y_x(\tau)\|_{L^\infty} \|p_{xx}(\tau)\|_{L^2}.
	$$
	We deduce
	\begin{equation}\label{eq:est18}
	\frac{d}{dt}\|p_x(\tau)\|^2 + \|p_{xx}(\tau)\|^2 \leq 2M^2 \|p_x(\tau)\|^2 + 2\|a_\tau(\tau)\|^2 \|y_{xx}(\tau)\|^2.
	\end{equation}
	Choose $\tau_1 > 2\epsilon$ and integrate the above inequality from $\tau_1$ to $\tau\in (\tau_1,T)$. We obtain
	$$
	\|p_x(\tau)\|^2 \leq \|p_x(\tau_1)\|^2 + 2M^2 \int_{2\varepsilon}^T\|p_x(\tau)\|^2\, d\tau + 2\int_{2\varepsilon}^T \|a_\tau(\tau)\|^2\, d\tau \esup_{\tau\in (2\varepsilon,T]}\|y_{xx}(\tau)\|^2.
	$$
	Now let $\tau\geq 3\varepsilon$ and integrate the above inequality over $\tau_1$ from $2\varepsilon$ to $3\varepsilon$. We obtain
	$$
	\varepsilon \|p_x(\tau)\|^2 \leq \int_{2\varepsilon}^T\|p_x(\tau_1)\|^2\, d\tau_1 + 2M^2 \varepsilon \int_{2\varepsilon}^T\|p_x(\tau)\|^2\, d\tau + 2\varepsilon \int_{2\varepsilon}^T \|a_\tau(\tau)\|^2\, d\tau \esup_{\tau\in (2\varepsilon,T]}\|y_{xx}(\tau)\|^2.
	$$
	It follows that $p_x = y_{\tau x} \in L^\infty(3\varepsilon,T;L^2)$, and, by \eqref{eq:est18}, $p_{xx} = y_{txx} \in L^2(3\varepsilon,T;L^2)$. It also follows from \eqref{eq:eq16} that $p_\tau = y_{\tau\tau} \in L^2(3\varepsilon,T;L^2)$.
	
	We have proved that $y_\tau \in L^2(3\varepsilon,T;H^2_{per})$ and $y_{\tau\tau} \in L^2(3\varepsilon,T;L^2_{per})$. This regularity implies that $y_\tau \in C([3\varepsilon,T];H^1_{per})$, whence $y_\tau \in C([0,1]\times [3\varepsilon,T])$. We have also proved that $y_{\tau x} \in L^2(3\varepsilon,T;H^1_{per})$, whence we deduce that  $y_{\tau x} \in L^2(3\varepsilon,T;C([0,1]))$ and $y_{\tau x}$ is $1$-periodic with respect to $x$ variable. It follows that $y_x\in C([0,1]\times [3\varepsilon,T])$ is also a $1$-periodic function with respect to $x$ variable. We deduce from \eqref{eq:adjoint_1_abs} that $y_{xx} \in C([0,1]\times [3\varepsilon,T])$. We have obtained enough smoothness of $y$ to use the strong maximum principle. 	
	
	Regularity \eqref{eq:Y_C} implies that $y(\tau)\in C([0,1])$ for every $\tau > 0$. By Lemma \ref{lemma:weak_max}
	 it follows that for every $\tau \in (0,T]$ there holds $\max_{x\in [0,1]}|y(x,\tau)| \leq \|y_0\|_{L^\infty}$. Without loss of generality we may assume that  $\max_{x\in [0,1]} y(x,\tau) \leq  \|y_0\|_{L^\infty}$ for every $\tau > 0$.
	 Assume that for some $\tau>0$ there holds $\max_{x\in [0,1]}|y(x,\tau)| = \|y_0\|_{L^\infty}$.
	 If there exists $x_0 \in (0,1)$ such that $y(x_0,\tau) = \|y_0\|_{L^\infty}$, then the strong maximum principle, see \cite[Theorem 2, page 168]{ProtterWeinberger} implies that $y(x,\tau) = \|y_0\|_{L^\infty}$ for every $(x,\tau) \in [0,1]\times [3\varepsilon,T]$. This means that
	$$
	y(3\varepsilon) - y(0) = \|y_0\|_{L^\infty} - y_0,
	$$
	and $\|y(3\varepsilon) - y(0)\| = \|\|y_0\|_{L^\infty} - y_0\| \neq 0$ as $y_0$ is not almost everywhere equal to a constant function. But $y\in C([0,T];L^2)$, so $\lim_{\varepsilon\to 0}\|y(3\varepsilon) - y(0)\| = 0$, a contradiction.
	We can hence exclude the case $y(x_0,\tau) = \|y_0\|_{L^\infty}$ for $x_0\in (0,1)$. Since $y$ is $1$-periodic with respect to $x$, we deduce that if $\max_{x\in [0,1]}|y(x,\tau)| = \|y_0\|_{L^\infty}$ then $y(0,\tau) = y(1,\tau) = 1$. We use   \cite[Theorem 3, page 170]{ProtterWeinberger} whence it follows that $y_x(0,\tau) < 0$ and $y_x(1,\tau) > 0$, a contradiction with $1$-periodicity of $y_x$ with respect to variable $x$. The proof is complete.
\end{proof}

 \begin{theorem}\label{theorem:l1_per}
 	Let $f\in L^\infty(\R;\dot{L}^2)$. The pullback attractor $\mathbb{A} = \{ A(t) \}_{t\in \R}$ established in Theorem 	\ref{thm:pullback_per} consists of a single eternal solution $A(t) = \{u(t)\}$ for every $t\in \R$. In other words there exists a unique eternal solution $u \in C_b(\R;\dot{H}^1)$ such that for every bounded set $B\in \mathcal{B}(\dot{L}^2)$ of initial data there holds
 	$$
 	\lim_{t_0 \to -\infty} \dist_{\dot{H}^1}(S(t,t_0)B,\{ u(t) \}) = 0.
 	$$
 \end{theorem}
 \begin{proof}
 	The first part of the argument follows the lines of the argument for the Dirichlet problem. Let $u,v: \R\to \dot{H}^1$ be two eternal solutions such that $\|u(t)\|_{\dot{H}^1} \leq M$ and $\|v(t)\|_{\dot{H}^1} \leq M$ for every $t\in \R$. Similar as in the Dirichlet case denote $a(t) = \frac{1}{2}(v(t) + v(t))$. Then $a\in C_b(\R,\dot{H}^1)$, $\|a(t)\|_{\dot{H}^1}\leq M$ for every $t\in \R$, $a_t\in L^2_{loc}(\R;L^2)$ and
 	$$
 	|a(x,t)| \leq \|a(t)\|_{L^\infty} \leq \|a_x(t)\| \leq M \quad \textrm{for every}\ \quad (x,t) \in [0,1]\times \R.
 	$$
 	Denote and $w(t) = u(t)-v(t)$. Then  $\|w(t)\|_{\dot{H}^1}\leq 2M$ for every $t\in \R$. Moreover, the function $w \in C_b(\R;\dot{H}^1)$ satisfies the equation
 	\begin{equation}\label{eq:forward_per}
 	w_t - w_{xx} + (a w)_x = 0,
 	\end{equation}
 	as well as the periodic conditions $w(0,t) = w(1,t)$, $w_x(0,1) = w_x(1,t)$ for every $t\in \R$.
 	Our aim is to prove that $u(t) = v(t)$, i.e., the only solution in $C_b(\R;\dot{H}^1)$ of the above equation such that $\|w(t)\|_{\dot{H}^1}\leq 2M$ for every $t$ is $w(t) = 0$. To this end we fix $t\in \R$ such that $w(t)$ is not identically zero, and consider \eqref{eq:forward_per} on the time interval $(t_0,t)$ for $t_0 < t$. As $w(t)$ is a continuous function of the variable $x$ we can define two open sets $A_+$ and $A_-$ by $A_+ = \{ x\in [0,1]\, :\ w(x,t) > 0 \}$, and $A_- = \{ x\in [0,1]\, :\ w(x,t) < 0 \}$ and  $z_0:[0,1] \to \R$ by
 	$$
 	z_0(x) = \chi_{A_+}(x) - \chi_{A_-}(x).
 	$$
 	It is clear that $z_0 \in L^2$ (but its mean is not necessarily equal to zero). Moreover
 	$$
 	\|w(t)\|_{L^1} = \int_0^1 z_0(x) w(x,t)\, dx.
 	$$
 	Similar as in the Dirichlet case we can define the adjoint problem which we solve backwards in time
 	\begin{align}
 	& z_t(s) + z_{xx}(s) + a(s)z_x(s) = 0\quad \textrm{for}\quad (x,s)\in (0,1)\times (t_0,t),\label{eq:backwards_per}\\
 	& z(0,s) = z(1,s) \quad \textrm{and} \quad z(0,s)_x= z(1,s)_x \quad \textrm{for}\quad s \in [t_0,t],\\
 	& z(t) = z_0.\label{eq:backwards3_per}
 	\end{align}
 	This problem has a unique weak solution $z \in L^2(t_0,t;H^1_{per})$ with $z_t\in L^2(t_0,t;(H^1_{per})')$. Note that $\dot{H}^1 \subset H^1_{per}$, the former being in fact a closed subspace of the latter. Proceeding exactly as in the Dirichlet case it follows that
 	\begin{equation}\label{eq:8_per}
 	(w(t_0),z(t_0)) = (w(t),z(t)) = (w(t),z_0) = \|w(t)\|_{\dot{L}^1}.
 	\end{equation}
 	Introducing the time $\tau = t-s$ the problem \eqref{eq:backwards_per}--\eqref{eq:backwards3_per} is equivalent to the following forward in time problem
 	\begin{align*}
 	&y_\tau(\tau) - y_{xx}(\tau) - a(t-\tau) y_x(\tau) = 0 \quad \textrm{for}\quad (x,\tau) \in (0,1)\times (0,t-t_0),\\
 	&y(0,\tau) = y(1,\tau)\quad \textrm{and}\quad  y_x(1,\tau) = y_x(0,\tau) \quad \textrm{for}\quad \tau \in (0,t-t_0),\\
 	& y(0) = z_0,
 	\end{align*}
 	namely, its solution $y$ is given by $y(\tau) = z(t-\tau)$ for $\tau \in [0,t-t_0]$. We will consider this problem on the whole positive semiaxis, which corresponds to taking arbitrarily small $t_0$, i.e
 	\begin{align}
 	&y_\tau(\tau) - y_{xx}(\tau) - a(t-\tau) y_x(\tau) = 0 \quad \textrm{for}\quad (x,\tau) \in (0,1)\times (0,\infty),\label{eq:adjoint_1}\\
 	&y(0,\tau) = y(1,\tau)\quad \textrm{and}\quad  y_x(1,\tau) = y_x(0,\tau) \quad \textrm{for}\quad \tau \in (0,\infty),\\
 	& y(0) = z_0,\label{eq:adjoint_3}
 	\end{align}
  It follows from \eqref{eq:8_per} and Lemma \ref{lemma:weak_max} that
 	 $$
 	 \|w(t)\|_{\dot{L}^1} \leq (w(t_0),z(t_0)) \leq \|w(t_0)\|_{\dot{L}^1} \|z(t_0)\|_{L^\infty} =  \|w(t_0)\|_{\dot{L}^1}  \|y(t-t_0)\|_{L^\infty} \leq \|w(t_0)\|_{\dot{L}^1}.
 	 $$
 	 Hence, $\|w(t)\|_{L^1}$ is a non-increasing function of $t$, whereas there exists numbers $0\leq c_1\leq c_2$ such that
 	 \begin{equation}\label{eq:convergence_t_infty}
 	 \lim_{t\to-\infty}\|w(t)\|_{L^1} = c_2\quad \textrm{and}\quad  \lim_{t\to\infty}\|w(t)\|_{L^1} = c_1
 	 \end{equation}
 	 As $w(t)$ is not identically zero, it follows that $c_2\neq 0$, we will show that this leads to a contradiction.
 	 Define $u^{\tau}(t) = u(t+\tau)$ and $v^\tau(t) = v(t+\tau)$. Then
 	\begin{equation}\label{eq:strong_u}
 	 	 u^\tau_t(t) - u^{\tau}_{xx}(t) + u^\tau(t)u_x^{\tau}(t) = f(t+\tau)\quad \textrm{for every}\quad t\in \R,
  	\end{equation}
 	 and
 	\begin{equation}\label{eq:strong_v}
 	 v^\tau_t(t) - v^{\tau}_{xx}(t) + v^\tau(t)v_x^{\tau}(t) = f(t+\tau) \quad \textrm{for every}\quad t\in \R,
  	\end{equation}
 	 that is, $u^{\tau}$ and $v^{\tau}$ are two eternal bounded in $\dot{H}^1$ solutions with the forcing terms $f(\cdot + \tau) = f^\tau$. Estimates of Section \ref{sec:periodic_estimates} imply that for every $t_1 < t_2$ there hold  the bounds
 	 $$
 	 \|u^\tau_x(t)\|+ \|v^\tau_x(t)\| \leq C,
 	 $$
 	 $$
 	 \int_{t_1}^{t_2} \|u^\tau_{xx}(t)\|^2\, dt + \int_{t_1}^{t_2} \|v^\tau_{xx}(t)\|^2\, dt \leq C,
 	 $$
 	 $$
 	 \int_{t_1}^{t_2} \|u^\tau_t(t)\|^2\, dt + \int_{t_1}^{t_2} \|v^\tau_t(t)\|^2\, dt \leq C,
 	 $$
 	 where the constant $C$ depends on $t_1-t_2$ but is independent of $\tau$. It is also clear that
 	 $$
 	 \|f^\tau(t)\| \leq C.
 	 $$
 	 We choose a sequence $\tau_n\to - \infty$. Using the diagonal argument it follows that there exists the subsequence of indexes, which we still denote by $n$, and functions $\overline{f}, \overline{v}, \overline{u}$ such that for every $t_1 < t_2$ the following convergences hold for $n\to \infty$
 	 \begin{align*}
 	 & f^{\tau_n} \to \overline{f} \quad \textrm{weakly}-*\ \textrm{in}\quad L^\infty(t_1,t_2;\dot{L}^2),\\
 	 & u^{\tau_n} \to \overline{u} \quad \textrm{weakly}-*\ \textrm{in}\quad L^\infty(t_1,t_2;\dot{H}^1)\ \ \textrm{and weakly in}\ \ L^2(t_1,t_2;\dot{H}^2),\\
 	 & v^{\tau_n} \to \overline{v} \quad \textrm{weakly}-*\ \textrm{in}\quad L^\infty(t_1,t_2;\dot{H}^1)\ \ \textrm{and weakly in}\ \ L^2(t_1,t_2;\dot{H}^2),\\
 	 & u^{\tau_n}_t \to \overline{u}_t \quad \textrm{and}\quad v^{\tau_n}_t \to \overline{v}_t \quad \textrm{weakly}\ \textrm{in}\quad L^2(t_1,t_2;\dot{L}^2),
 	 \end{align*}
 	 and
 	  $$
 	 u^{\tau_n}(t) \to \overline{u}(t)\quad \textrm{and}\quad v^{\tau_n}(t) \to \overline{v}(t) \quad \textrm{weakly in}\quad \dot{H}^1\quad \textrm{and strongly in}\quad \dot{L}^2\quad \textrm{for every}\quad t\in \R.
 	 $$
 	 These convergences allow us to pass to the limit with $\tau_n$ to minus infinity in \eqref{eq:strong_u} and \eqref{eq:strong_v}, whence it follows that
 	 $$
 	 \overline{u}, \overline{v} \in C_b(\R;\dot{H}^1)\cap L^2_{loc}(\R;\dot{H}^2), \overline{u}_t, \overline{v}_t \in L^2_{loc}(\R;\dot{L}^2)
 	 $$
 	 satisfy (in strong sense) the equations
 	 \begin{align*}
 	 & \overline{u}_t(t) -\overline{u}_{xx}(t) + \overline{u}(t)\overline{u}_x(t) = \overline{f}(t),\\
 	 & \overline{v}_t(t)-\overline{v}_{xx}(t) + \overline{v}(t)\overline{v}_x(t) = \overline{f}(t).
 	 \end{align*}
 	 Define $\overline{w} = \overline{u}-\overline{v}$ and $\overline{a} = \frac{1}{2}(\overline{u}+\overline{v})$. Then $\|\overline{a}(t)\|_{\dot{H}^1} \leq M$ and $\|\overline{w}(t)\|_{\dot{H}^1} \leq 2M$ for every $t \in \R$  and we have the following regularities and convergences
 	 $$
 	 \overline{w},\overline{a} \in C_b(\R;\dot{H}^1)\cap L^2_{loc}(\R;\dot{H}^2), \overline{w}_t, \overline{a}_t \in L^2_{loc}(\R;\dot{L}^2),
 	 $$
 	 $$
 	 w(\dot + \tau_n) = w^{\tau_n} \to \overline{w}  \quad \textrm{weakly}-*\ \textrm{in}\quad L^\infty(t_1,t_2;\dot{H}^1)\ \ \textrm{and weakly in}\ \ L^2(t_1,t_2;\dot{H}^2),
 	 $$
 	 $$
 	 w^{\tau_n}(t) \to \overline{w}(t) \quad \textrm{weakly in}\quad \dot{H}^1\quad \textrm{and strongly in}\quad \dot{L}^2\quad \textrm{for every}\quad t\in \R.
 	 $$
 	 Moreover, $\overline{w}$ satisfies is strong sense the equation
 	 $$
 	 \overline{w}_t(t) - \overline{w}_{xx}(t) + (\overline{a}(t)\overline{w}(t))_x = 0.
 	 $$
 	 The convergence \eqref{eq:convergence_t_infty} implies that
 	 $$
 	 \|w^{\tau_n}(t)\|_{L^1}\to c_2 \quad \textrm{for every}\quad t\in \R,
 	 $$
 	 and hence
 	 $$
 	 \|\overline{w}(t)\|_{L^1} = c_2 \quad \textrm{for every}\quad t\in \R.
 	 $$
 	 Choose real numbers $t_0 < t$ and define $z_0(x) = \chi_{\{x\in [0,1]\, :\ \overline{w}(x,t)>0\}} - \chi_{\{x\in [0,1]\, :\ \overline{w}(x,t)<0\}}$.
 	Analogously to \eqref{eq:adjoint_1}--\eqref{eq:adjoint_3} we formulate the adjoint problem
 	 	\begin{align}
 	&y_\tau(\tau) - y_{xx}(\tau) - \overline{a}(t-\tau) y_x(\tau) = 0 \quad \textrm{for}\quad (x,\tau) \in (0,1)\times (0,\infty),\label{eq:adjoint_1__norm_c}\\
 	&y(0,\tau) = y(1,\tau)\quad \textrm{and}\quad  y_x(1,\tau) = y_x(0,\tau) \quad \textrm{for}\quad \tau \in (0,\infty),\\
 	& y(0) = z_0.\label{eq:adjoint_3__norm_c}
 	\end{align}
 	As in \eqref{eq:8_per}, as $\overline{w}(t)$ is not identically zero, we get $\|z_0\|_{L^\infty} = 1$, and
 	$$
 	c_2  = \|\overline{w}(t)\|_{\dot{L}_1} = (y(t-t_0),w(t)) \leq \|y(t-t_0)\|_{L^\infty}\|w(t_0)\|_{L^1} = c_2\|y(t-t_0)\|_{L^\infty}.
 	$$
 	As $c_2\neq 0$ and $t_0$ is arbitrary it follows that
 	$$
 	\|y(\tau)\|_{L^\infty} \geq 1\quad \textrm{for every} \quad \tau \geq 0.
 	$$
 	But Lemma \ref{lem:strong_max} leads us to the conclusion that
 	$$
 	\|y(\tau)\|_{L^\infty} < 1\quad \textrm{for every} \quad \tau > 0,
 	$$
 	which is a contradiction and the proof is complete
 \end{proof}

We can repeat the argument of the above theorem passing with $\tau_n$ to $+\infty$ instead of $-\infty$ to deduce that the unique complete trajectory also attracts all trajectories in future. That is, we get the following result
 \begin{theorem}\label{theorem:l1_per_2}
	Let $f\in L^\infty(\R;\dot{L}^2)$.  The unique eternal strong solution $u \in C_b(\R;\dot{H}^1)$ satisfies
	$$
	\lim_{t\to\infty} \|S(t,t_0)u_0 - u(t)\| = 0
	$$
	for every $t_0 \in \R$ and $u_0 \in \dot{L}^2$.
	
\end{theorem}

\begin{remark}
	Similar as in the Dirichlet case it is clear that the $T$-periodicity of the non-autonomous forcing $f \in L^\infty(\R;\dot{L}^2)$ implies that the unique eternal solution bounded in $\dot{H}^1$ which attracts all trajectories in pullback and forward sense is also $T$-periodic.
\end{remark}

\end{document}